\documentclass[12pt]{amsart}
\usepackage[leqno]{amsmath}
\usepackage[psamsfonts]{amssymb}
\usepackage{amsthm}
\usepackage[mathscr,psamsfonts]{eucal}
\usepackage{cmmib57}
\usepackage{amscd}
\usepackage{url}
\usepackage{color}

\IfFileExists{Ueur.fd}{\input{Ueur.fd}}{\input{ueur.fd}}
\IfFileExists{Ueur57.fd}{\input{Ueur57.fd}}{\input{ueur57.fd}}
\DeclareMathAlphabet{\eurm}{U}{eur}{m}{n}
\DeclareMathAlphabet{\eubf}{U}{eur}{b}{n}
\font\cyr=wncyr10 scaled 1200

\newtheoremstyle
{MyThm}
{10pt}
{10pt}
{\itshape}
{\parindent}
{\bfseries}
{.}
{.5em}
{}
\swapnumbers
\theoremstyle{MyThm}
\newcounter{assump}
\newtheorem{Assumption}{Assumption}[assump]

\newcounter{postul}
\newtheorem{Postulate}{Postulate}[postul]

\newtheorem{Caution}{Caution}[section]
\newtheorem{Convention}[Caution]{Convention}
\newtheorem{Corollary}[Caution]{Corollary}
\newtheorem{Definition}[Caution]{Definition}
\newtheorem{Example}[Caution]{Example}
\newtheorem{Exercise}[Caution]{Exercise}
\newtheorem{Lemma}[Caution]{Lemma}
\newtheorem{Notation}[Caution]{Notation}
\newtheorem{Note}[Caution]{Note}
\newtheorem{Problem}[Caution]{Problem}
\newtheorem{Proposition}[Caution]{Proposition}
\newtheorem{Remark}[Caution]{Remark}
\newtheorem{Theorem}[Caution]{Theorem}
\newcommand{\bAs}{\begin{Assumption}\em}
\newcommand{\eAs}{\end{Assumption}}
\newcommand{\bCa}{\begin{Caution}\em}
\newcommand{\eCa}{\end{Caution}}
\newcommand{\bCr}{\begin{Corollary}\em}
\newcommand{\eCr}{\end{Corollary}}
\newcommand{\bCv}{\begin{Convention}\em}
\newcommand{\eCv}{\end{Convention}}
\newcommand{\bDf}{\begin{Definition}\em}
\newcommand{\eDf}{\end{Definition}}
\newcommand{\bDr}{\begin{Exercise}\em}
\newcommand{\eDr}{\end{Exercise}}
\newcommand{\bEx}{\begin{Example}\em}
\newcommand{\eEx}{\end{Example}}
\newcommand{\bLm}{\begin{Lemma}\em}
\newcommand{\eLm}{\end{Lemma}}
\newcommand{\bNo}{\begin{Notation}\em}
\newcommand{\eNo}{\end{Notation}}
\newcommand{\bNt}{\begin{Note}\em}
\newcommand{\eNt}{\end{Note}}

\newcommand{\bPb}{\begin{Problem}\em}
\newcommand{\ePb}{\end{Problem}}
\newcommand{\bPf}{\begin{proof}}
\newcommand{\ePf}{\end{proof}}
\newcommand{\bpf}{\bPf}
\newcommand{\epf}{\ePf}
\newcommand{\bPr}{\begin{Proposition}\em}
\newcommand{\ePr}{\end{Proposition}}
\newcommand{\bPs}{\begin{Postulate}\em}
\newcommand{\ePs}{\end{Postulate}}
\newcommand{\bRm}{\begin{Remark}\em}
\newcommand{\eRm}{\end{Remark}}

\newcommand{\bTh}{\begin{Theorem}}
\newcommand{\eTh}{\end{Theorem}}
\newcommand{\bEq}{\begin{equation}}
\newcommand{\eEq}{\end{equation}}
\newcommand{\beq}{\begin{equation*}}
\newcommand{\eeq}{\end{equation*}}

\newcommand{\bal}{\begin{align*}}
\newcommand{\bAl}{\begin{align}}
\newcommand{\bat}{\begin{alignat*}}
\newcommand{\bAt}{\begin{alignat}}
\newcommand{\bml}{\begin{multline*}}
\newcommand{\bMl}{\begin{multline}}
\newcommand{\bgt}{\begin{gather*}}
\newcommand{\bGt}{\begin{gather}}
\newcommand{\bCd}{\bEq\begin{CD}}
\newcommand{\eCd}{\end{CD}\eEq}
\newcommand{\bcd}{\beq\begin{CD}}
\newcommand{\ecd}{\end{CD}\eeq}
\newcommand{\bdg}{\beq\begin{diagram}}
\newcommand{\edg}{\end{diagram}\eeq}
\newcommand{\bDg}{\bEq\begin{diagram}}
\newcommand{\eDg}{\end{diagram}\eEq}
\newcommand{\bmt}{\left(\begin{matrix}}
\newcommand{\emt}{\end{matrix}\right)}
\newcommand{\bcn}{\begin{center}}
\newcommand{\ecn}{\end{center}}
\newcommand{\ben}{\begin{enumerate}}
\newcommand{\een}{\end{enumerate}}
\newcommand{\btb}{\begin{tabbing}}
\newcommand{\etb}{\end{tabbing}}
\newcommand{\bfz}{\begin{footnotesize}}
\newcommand{\efz}{\end{footnotesize}}
\newcommand{\bsz}{\begin{scriptsize}}
\newcommand{\esz}{\end{scriptsize}}

\newcommand{\bsb}
{\vspace{-0.8cm}
\begin{alignat*}{2}
& \qquad\qquad\qquad\qquad\qquad\qquad\qquad\qquad\qquad\qquad
&&\qquad\qquad\qquad\qquad\qquad\qquad\qquad\qquad\qquad
\\}
\newcommand{\Rn}{{\B R}}

\newcommand{\Al}{\forall}

\newcommand{\h}{\hbar}
\newcommand{\1}{\mathbf 1}

\newcommand{\der}{\partial}

\newcommand{\nab}{\nabla}

\newcommand{\Fla}{^{\flat}{}}
\newcommand{\Sha}{^{\sharp}{}}
\newcommand{\ShaN}{^{\sharp_{\f N}}{}}

\newcommand{\lang}{\langle}
\newcommand{\rang}{\rangle}

\newcommand{\mto}{\mapsto}
\newcommand{\hto}{\hookrightarrow}
\newcommand{\sub}{\subset}
\newcommand{\com}{\circ}

\newcommand{\car}{\times}
\newcommand{\ten}{\otimes}
\newcommand{\drs}{\oplus}

\newcommand{\wed}{\wedge}

\DeclareMathOperator{\con}{\lrcorner}

\newcommand{\eqv}{\,\equiv\,}
\newcommand{\seq}{\,\simeq\,}
\newcommand{\nid}{\not\equiv}
\DeclareMathOperator{\byd}{\,{\rm =}{\raisebox{.092ex}{\rm :}}\,}

\newcommand{\ucar}[1]{\underset{#1}{\times}}

\newcommand{\fr}[2]{\frac{#1}{#2}\,}
\newcommand{\tfr}[2]{\tfrac{#1}{#2}\,}

\newcommand{\col}[3]{_{#1}{}^{#2}{}_{#3}}

\newcommand{\Ga}[2]{_{#1}{}^{#2}_0}

\newcommand{\Gaa}[3]{_{#1}{}^{#2}_0{}^0_{#3}}

\newcommand{\ga}[1]{_0{}^{#1}_0}

\newcommand{\ENDE}{{\,\text{\footnotesize\qedsymbol}}}

\newcommand{\ssep}[1]{{\qquad\text{\rm{#1}}\qquad}}
\newcommand{\st}{\;|\;}

\newcommand{\ar}[1]{\url{http://arXiv.org/abs/#1}}

\newcommand{\im}{{{}{\rm im \, }}}

\DeclareMathOperator{\Grass}{{{Grass}}}

\DeclareMathOperator{\Span}{{{span}}}

\DeclareMathOperator{\fib}{{{fib}}}

\DeclareMathOperator{\map}{{{map}}}

\newcommand{\f}[1]{{\boldsymbol{#1}}}

\newcommand{\ul}[1]{{\underline{#1}}}
\newcommand{\ba}[1]{{{\bar{#1}}}}

\newcommand{\ch}[1]{{\check{#1}}}
\newcommand{\wch}[1]{{\overset{\vee}{#1}}}

\newcommand{\ti}[1]{{\tilde{#1}}}
\newcommand{\wti}[1]{{\widetilde{#1}}}
\newcommand{\dt}[1]{{\dot{#1}}}

\newcommand{\br}[1]{\breve{#1}{}}

\newcommand{\ltag}[1]{\leqno{\quad\text{#1)}}}
\newcommand{\bma}{\left(\begin{matrix}}
\newcommand{\ema}{\end{matrix}\right)}

\newcommand{\M}[1]{{\mathscr{#1}}}

\newcommand{\B}[1]{{\mathbb{#1}}}
\newcommand{\baB}[1]{{\bar{{\mathbb{#1}}}}}

\newcommand{\K}[1]{\text{\cyr{#1}}}

\newcommand{\alp}{\alpha}
\newcommand{\bet}{\beta}
\newcommand{\gam}{\gamma}
\newcommand{\del}{\delta}

\newcommand{\tht}{\theta}

\newcommand{\lam}{\lambda}

\newcommand{\sig}{\sigma}

\newcommand{\ome}{\omega}
\newcommand{\Gam}{\Gamma}

\newcommand{\Lam}{\Lambda}

\newcommand{\Ome}{\Omega}

\allowdisplaybreaks[4]
\begin{document}

\title[Generalized geometrical structures]
{Generalized geometrical structures
\\
of odd dimensional manifolds}

\author[J. Jany\v{s}ka, M. Modugno]
	{Josef Jany\v{s}ka and Marco Modugno}
\address{
{\ }
\newline
Department of Mathematics and Statistics, Masaryk University
\newline
Jan\'a\v{c}kovo n\'am 2a, 602 00 Brno, Czech Republic
\newline
E-mail: {\tt janyska@math.muni.cz}
\newline
{\ }
\newline
Department of Applied Mathematics, Florence University
\newline
Via S. Marta 3, 50139 Florence, Italy
\newline
E-mail: {\tt marco.modugno@unifi.it}
}



\thanks{
This research has been supported
by the Ministry of Education of the Czech Republic under the project
MSM0021622409,
by the Grant agency of the Czech Republic under the project GA
201/05/0523,
by MIUR of Italy under the project PRIN 2005 ``Simmetrie e
Supersimmetrie Classiche e Quantistiche",
by GNFM of INdAM
and by Florence University.
}
\begin{abstract}
 We define an almost--cosymplectic--contact structure which generalizes
co\-symplectic and contact structures of an odd dimensional manifold.
 Analogously, we define an almost--coPoisson--Jacobi
structure which  generalizes  a Jacobi structure.
 Moreover, we study relations between these structures and
analyse the associated algebras of functions.

 As examples of the above structures, we present geometrical dynamical
structures of the phase space of a general relativistic particle,
regarded as the 1st jet space of motions in a spacetime.
 We describe geometric conditions by which a metric and a connection
of the phase space yield cosymplectic and dual coPoisson structures,
in case of a spacetime with  absolute time (a Galilei spacetime), or
almost--cosymplectic--contact and dual  almost--coPoisson--Jacobi
structures, in case of a spacetime without absolute time (an
Einstein spacetime).
\end{abstract}
\maketitle
\section*{Introduction}
 In \cite{JadJanMod98, JanMod95, JanMod99, JanMod02, JanMod06}
we studied geometrical structures on the phase space of a spacetime
naturally induced (in the sense of \cite{KolMicSlo93}) by a metric
and  a phase connection.
 Some of these structures are well known and some are less standard.
 In the present paper, we generalize these structures on odd
dimensional manifolds and study general properties of such
structures. 

 First, in Section 1, we recall some standard structures and
introduce new structures, namely almost--cosymplectic--contact,
coPoisson and  almost--coPoisson--Jacobi structures.
 In Section 2 we study algebras of functions which are associated
with the new geometrical structures.

 As examples of the above new structures, we study the geometrical
structures on the phase space of a spacetime.
 Actually, the geometric objects arising in Section \ref{Galilei
spacetime}, in the framework of the Galilei's phase space 
\cite{JadJanMod98, JanMod99, JanMod02},
involve mainly the concepts of cosymplectic and (regular) 
coPoisson structures.
 On the other hand, the analogous geometric objects arising in
Section 
\ref{Einstein spacetime},
in the framework of the Einstein's phase space \cite{JanMod95,
JanMod06}, involve mainly the concepts of
almost--cosymplectic--contact and  almost--coPoisson--Jacobi
structures (eventually contact and Jacobi structures).
\section{Geometrical structures}
\setcounter{equation}{0}
 We use the inner product 
$i$ 
of 
$k$--vectors with 
$r$--forms defined by
$i_{X_1 \wed \dots X_k} \bet = i_{X_k} \dots i_{X_1} \bet \,,$
for each 
$r$-form 
$\bet$ 
and 
$k$
vector fields 
$X_1, \dots, X_k \,,$
with
$k \le r \,.$
 We use the same symbol for the dual inner product of $k$-forms with
$r$-vectors.

 For the Schouten bracket we use the identity, 
\cite{LibMar87, Lich78, Vai94},
\beq
i_{[P,Q]}\bet = 
(-1)^{q(p+1)} i_Pdi_Q\bet + (-1)^p i_Qdi_P\bet - i_{P\wed Q}d\bet \,,
\eeq
for each $p$--vector 
$P \,,$ 
$q$--vector 
$Q$ 
and 
$(p+q-1)$--form 
$\bet \,.$
 In particular, for each vector field 
$E$ 
and 2--vector 
$\Lam \,,$
we have
$i_{[E,\Lam]} \bet = i_E d i_\Lam \bet - i_\Lam d i_E \bet \,,$
for each closed 2-form $\bet \,,$ 
and
$i_{[\Lam,\Lam]}\bet = 2\, i_\Lam d i_\Lam \bet \,,$
for each closed 3--form 
$\bet \,.$

\smallskip

 In what follows,
$\f M$
is a 
$(2n + 1)$--dimensional smooth manifold.
\subsection{Covariant and contravariant pairs}
\label{Contravariant and covariant pairs}
\bDf
 We define a \emph{covariant pair} to be a pair
$(\ome, \Ome)$
consisting of a 1--form
$\ome$
and
a 2--form
$\Ome$
of constant rank
$2r \,,$
with
$0 \le r \le n \,,$
such that 
$\ome \wed \Ome^r \nid 0 \,,$ 
and a \emph{contravariant pair} to be a pair
$(E, \Lam)$
consisting of a vector field
$E$
and a 2--vector
$\Lam$
of constant rank
$2s \,,$
with
$0 \le s \le n \,,$
such that 
$E \wed \Lam^s \nid 0 \,.$
 Thus, by definition, we have
$\Ome^r \nid 0 \,,\; \Ome^{r+1} \eqv 0$
and
$\Lam^s \nid 0 \,,\; \Lam^{s+1} \eqv 0 \,.$

 We say that the pairs
$(\ome, \Ome)$
and
$(E, \Lam)$
are \emph{regular} if, respectively,
\beq
\ome \wed \Ome^n \nid 0
\ssep{and}
E \wed \Lam^n \nid 0 \,.
\eeq
\vglue-1.5\baselineskip{ }\hfill\ENDE
\eDf

 Let us consider a covariant pair
$(\ome,\Ome)$
and a contravariant pair
$(E,\Lam) \,.$

 We define the following linear maps and subspaces
\bat{3}
&\Ome\Fla : T\f M \to T^*\f M : X \mto X\Fla \byd i_X \Ome \,,
&&\qquad
&&\Lam\Sha : T^*\f M \to T\f M : \alp \mto \alp\Sha \byd i_\alp \Lam
\,,
\\
&\lang\ome\rang \byd \{\lam \ome \st \lam \in \Rn\} \sub T^*\f M \,,
&&\qquad
&&\lang E \rang \byd \{\lam E \st \lam \in \Rn\} \sub T\f M \,,
\\
&\ker E \byd \{\alp \in T^*\f M \st \alp(E) = 0\} \,,
&&\qquad
&&\ker \ome \byd \{X \in T\f M \st \ome(X) = 0\} \,.
\end{alignat*}

 We have
$\dim \, (\im\Ome\Fla) = 2r$
and
$\dim \, (\im\Lam\Sha) = 2s \,.$

  If
$(\ome,\Ome)$
is regular, then
$r = n \,,$
$\dim \, (\im\Ome\Fla)  = 2n \,,$
$\dim \, (\ker\Ome\Fla) = 1 \,,$
$\dim \, (\ker\ome) = 2n \,.$

  If
$(E,\Lam)$
is regular, then
$s = n \,,$
$\dim \, (\im\Lam\Sha)  = 2n \,,$
$\dim \, (\ker\Lam\Sha) = 1 \,,$
$\dim \, (\ker E) = 2n \,.$
\subsection{Structures given by covariant pairs}
\label{Structures given by covariant pairs}
 According to \cite{Lich78}, a  
\emph{pre cosymplectic structure} on 
$\f M$ 
is defined by a regular covariant pair
$(\ome, \Ome) \,.$

 Two distinguished types of pre cosymplectic  structures appear in
the literature.
 Namely, we recall that a \emph{cosymplectic structure}
\cite{deLTuy96}
and a \emph{contact structure} 
\cite{LibMar87} are defined by a covariant pair
$(\ome,\Ome)$ 
such that, respectively,
\bGt
d\ome = 0 \,,
\qquad
d\Ome = 0 \,,
\qquad
\ome \wed \Ome^n \nid 0 \,,
\\
\Ome = d\ome \,,
\qquad
\ome \wed \Ome^n \nid 0 \,.
\end{gather}

 Thus, a contact structure is characterised just by a 1--form 
$\ome$
such that
\beq
\ome \wed (d\ome)^n \nid 0 \,.
\eeq

 We can easily generalize the above structures in the following way.

\bDf
 We define an \emph{almost--cosymplectic--contact structure} to be a
covariant pair 
$(\ome,\Ome)$
such that
\beq
d\Ome = 0 \,,
\qquad
\ome \wed \Ome^n \nid 0 \,. \eqno{\ENDE}
\eeq
\eDf 

 Clearly, for 
$d\ome = 0$ 
we obtain a cosymplectic structure and for
$\Ome = d\ome$ 
a contact structure. 
 So, almost--cosymplectic--contact structures are regular structures
which generalize both cosymplectic and contact structures.
\subsection{Structures given by contravariant pairs}
\label{Structures given by contravariant pairs}
 Two distinguished types of contravariant pairs appear in the
literature.

 Namely, we recall that a \emph{Jacobi structure} is defined by a
contravariant pair
$(E, \Lam)$
such that
\beq
[E, \Lam] = 0 \,,
\qquad
[\Lam, \Lam] = - 2 E \wed \Lam \,,
\eeq
where
$[\,,]$
denotes the Schouten bracket.

 In the particular case when 
$E = 0 \,,$
we obtain
\beq
[\Lam, \Lam] = 0
\eeq
and the pair
$(E, \Lam) \byd (0, \Lam)$  
is called \emph{Poisson structure}.

 On the other hand, in the particular case when
$\Lam = 0 \,,$
we obtain
$[E,\Lam] = 0$
and 
$[\Lam,\Lam] = 0$
and the pair
$(E,\Lam) \byd (E, 0)$ 
is called \emph{trivial structure}.

 In the following we assume 
$E \nid 0$ 
and 
$\Lam \nid 0 \,.$

\bRm
 In the literature (see for instance \cite{Lich78}) the condition 
$E \wed \Lam^s \nid 0$
is considered just as a possible non necessary property of the Jacobi
pair
$(E, \Lam) \,.$
 So, our definition is a little more restrictive; however, the
assumption
$E \wed \Lam^s \nid 0$
is quite reasonable and it is needed for our subsequent developments. 

 In the literature (see for instance
\cite{LibMar87, Lich78, Vai94}) 
the Jacobi structure is usually defined by the identities 
$[E, \Lam] = 0 \,,\;\; [\Lam, \Lam] =  2 E \wed \Lam \,. $
 The difference in the sign in the second identity, with respect to
our definition, is caused by the different convention on the inner
product, hence by the different sign in definition of 
$\Lam\Sha \,.$
\hfill\ENDE
\eRm

 In order to exhibit a certain symmetry between geometric structures
given by covariant and contravariant pairs, we introduce the following
notions.

\smallskip

\bDf
 We define a \emph{pre coPoisson structure} to be a contravariant
pair
$(E, \Lam) \,.$

 In particular,  a \emph{coPoisson structure} is defined by a
contravariant pair
$(E, \Lam)$
such that
$$
[E, \Lam] = 0 \,,
\qquad
[\Lam, \Lam] = 0 \,.
\eqno{\ENDE}
$$
\eDf

\bDf\label{Df1.5}
 We define an \emph{almost--coPoisson--Jacobi structure} to be a
3--plet\linebreak
$(E, \, \Lam, \, \ome) \,,$  
where
$(E,\Lam)$  
is a contravariant pair and
$\ome$
a 1--form, such that  
\beq
[E,\Lam] = - E \wed \Lam\Sha(L_E\ome) \,,
\quad
[\Lam,\Lam] = 2 \, E \wed (\Lam\Sha \ten \Lam\Sha)(d\ome) \,,
\qquad
i_E \ome = 1\,,
\quad
i_\ome\Lam = 0 \,.
\eeq

 The 1-form 
$\ome$ 
is said to be the \emph{fundamental 1-form} of the
almost--coPoisson--Jacobi structure.
\hfill\ENDE
\eDf

\bRm
 Almost--coPoisson--Jacobi structures generalize both coPoisson
and Jacobi structures.

 Indeed, if
$d\ome = 0 \,,$
then we have
$L_E\ome = i_E d\ome = 0 \,,$ 
hence from Definition \ref{Df1.5} we obtain
$[E,\Lam] = 0$ 
and 
$[\Lam,\Lam] = 0 \,,$
i.e. 
$(E,\Lam)$ 
turns out to be a coPoisson structure.

 Moreover, if 
$L_E\ome = 0$
and
$(\Lam\Sha\ten\Lam\Sha)(d\ome) = - \Lam \,,$
then we obtain
$[E, \Lam] = 0$
and
$[\Lam, \Lam] = - 2 E \wed \Lam \,,$
i.e.
$(E, \Lam)$
turns out to be a Jacobi structure.
\hfill\ENDE
\eRm

\bPr
 Let
$(E,\Lam)$
be a regular contravariant pair.
 Then, there exists a unique 1--form 
$\ome \,,$
such that 
$i_\ome (E\wed \Lam^n) = \Lam^n \,.$
Indeed, such an
$\ome$ 
satisfies the equalities
$i_E\ome = 1$ 
and 
$i_\ome\Lam = 0 \,.$

 Thus, the 3--plet
$(E,\Lam, \ome)$
turns out to be an almost--coPoisson--Jacobi structure if and only if
$\,[E,\Lam] = - E \wed \Lam\Sha(L_E\ome)\,$
and
$\,[\Lam,\Lam] = 2 \, E \wed (\Lam\Sha \ten \Lam\Sha)(d\ome) \,.$
\hfill\ENDE
\ePr

 Thus, a regular almost--coPoisson--Jacobi structure can be defined just as a
suitable contravariant pair
$(E, \Lam) \,,$
as the additional 1--form
$\ome$
is naturally determined by the above pair itself.
\subsection{Dual structures}
\label{Dual structures}
 Let us consider a covariant pair
$(\ome,\Ome)$
and a contravariant pair
$(E,\Lam) \,.$

\bDf
 The pairs
$(\ome,\Ome)$
and
$(E,\Lam)$
are said to be
\emph{mutually dual} if they are regular, the maps
\beq
\Ome\Fla_{|\im (\Lam\Sha)} : \im (\Lam\Sha) \to
\im (\Ome\Fla) \sub T^*\f M
\ssep{and}
\Lam\Sha_{|\im (\Ome\Fla)} : \im (\Ome\Fla) \to
\im (\Lam\Sha) \sub T\f M
\eeq
are isomorphisms and
\beq
(\Ome\Fla_{|\im (\Lam\Sha)})^{-1} = 
\Lam\Sha_{|\im (\Ome\Fla)} \,,
\quad
(\Lam\Sha_{|\im (\Ome\Fla)})^{-1} = 
\Ome\Fla_{|\im (\Lam\Sha)} \,,
\quad 
i_E\Ome = 0 \,,
\quad 
i_\ome\Lam = 0 \,,
\quad 
i_E\ome = 1 \,.\hfill\ENDE
\eeq
\eDf

\bTh\label{Th-dual pairs}\cite{Lich78}
 The relation of duality yields a bijection between regular covariant
pairs
$(\ome,\Ome)$ 
and regular contravariant pairs
$(E,\Lam) \,.$
\hfill\ENDE
\eTh

 Thus, the geometric structures given by dual covariant and
contravariant pairs are essentially the same.

 In the literature 
$E$ 
is called the \emph{fundamental vector field} \cite{Lich78}, or the 
\emph{Reeb vector field} \cite{Ree52},
and 
$\Lam$ 
the \emph{fundamental 2-tensor} of
$(\ome,\Ome) \,.$

\bigskip

 Now, let us consider dual pairs
$(\ome,\Ome)$ 
and 
$(E,\Lam)$
and state some results.

\bLm
 We have
\beq
\lang E \rang = \ker \Ome\Fla \,,
\quad
\im (\Lam\Sha) = \ker \ome
\ssep{and}
\lang\ome\rang = \ker \Lam\Sha \,,
\quad
\im (\Ome\Fla) = \ker E \,.
\eeq
\eLm

\bpf
1) We have
$\lang E\rang \sub \ker \Ome\Fla \,;$
hence,
$\dim \, (\ker \Ome\Fla) = 1 = \dim \, \lang E \rang$
implies
$\lang E\rang = \ker \Ome\Fla \,.$

 If
$X \in \sec(\f M, \im (\Lam\Sha)) \,,$
then there exists
$\alp \in \sec(\f M, T^*\f M) \,,$
such that
$i_\alp \Lam = X \,;$
hence,
\beq
\ome(X) = \ome (i_\alp \Lam) = \Lam (\alp, \ome) = 
- i_\alp \Lam\Sha(\ome) = 0 \,,
\ssep{hence}
X \in \sec(\f M, \ker \ome) \,.
\eeq

 Then, 
$\dim \, (\im\Lam\Sha) = 2n = \dim \, (\ker \ome)$
implies
$\im (\Lam\Sha) = \ker \ome \,.$

2) In the same way we prove the other two identities.
\epf

\bPr
 We have the splittings
\beq
T\f M = \lang E \rang \drs \im (\Lam\Sha)
\ssep{and}
T^*\f M = \lang\ome\rang \drs \im (\Ome\Fla) \,.
\eeq

 Accordingly, for each
$X \in \sec(\f M, T\f M)$
and
$\alp \in \sec(\f M, T^*\f M) \,,$
we have the splittings
\beq
X = \ome(X) \, E + \big(X - \ome(X) \, E\big)
\ssep{and}
\alp = \alp(E) \, \ome + \big(\alp - \alp(E) \, \ome\big) \,,
\eeq

 Thus, the maps
\beq
\Lam\Sha \com \Ome\Fla : T\f M \to \im (\Lam\Sha)
\ssep{and}
\Ome\Fla \com \Lam\Sha : T^*\f M\to \im (\Ome\Fla)
\eeq
are the ``orthogonal" projections of the splittings of
$T\f M$
and 
$T^*\f M \,.$
\ePr

\bpf
 The equalities
$\dim \lang E\rang + \dim \im (\Lam\Sha) = 1 + 2 n$
and 
$\lang E\rang \cap \im (\Lam\Sha) = \lang E\rang \cap \ker\ome = 0$
yield
$T\f M = \lang E\rang \drs \im (\Lam\Sha) \,.$

 Clearly, we have
\beq
\ome(X) \, E \in \sec(\f M, \lang E \rang) \,,
\quad
X - \ome(X) \, E \in \sec(\f M, \im (\Lam\Sha)) = 
\sec(\f M, \ker \ome) \,.
\eeq

 Then, we obtain
\beq
X - \ome(X) \, E = (\Lam\Sha \com \Ome\Fla) (X - \ome(X) \, E) =
(\Lam\Sha \com \Ome\Fla) (X) \,.
\eeq

 The dual result can be obtained in the same way.  
\epf

\bPr
 For each 
$X, Y \in \sec(\f M, T\f M)$ 
and
$\alp, \bet \in \sec(\f M, T^*\f M) \,,$ 
we have
\bEq
\Ome(\alp\Sha,\bet\Sha) = - \Lam(\alp,\bet)
\ssep{and}
\quad \Lam(X\Fla,Y\Fla) = - \Ome(X,Y)\,,
\eEq
i.e.
\bEq
(\Lam\Sha\ten\Lam\Sha)(\Ome) = - \Lam
\ssep{and}
\quad (\Ome\Fla\ten\Ome\Fla)(\Lam) = - \Ome\,.
\eEq
\ePr

\bpf
 We have
\beq
\Ome\big(\Lam\Sha(\alp),\Lam\Sha(\bet)\big) = 
i_{\Lam\Sha(\bet)} \Ome\Fla\big(\Lam\Sha(\alp)\big) =
i_{\Lam\Sha(\bet)} \big(\alp - \alp(E) \, \ome\big) =
\Lam\big(\bet \,,\, \alp - \alp(E)\, \ome\big) =
- \Lam(\alp,\bet) \,.
\eeq

 The second identity can be proved in the same way. 
\epf

\bLm\label{Lemma: expression of d Omega}
 Let us consider the functions
$f, g, h \in \map(\f M, \Rn) \,,$
the \emph{closed} forms
$\alp, \bet, \gam \in \sec(\f M, T^*\f M) \,,$
and the induced vector fields
$X,Y,Z \in \sec(\f M, T\f M) \,,$
given by
\bEq\label{splitting of X, Y, Z}
X \byd \alp\Sha + f E \,,
\qquad
Y \byd \bet\Sha + g E \,,
\qquad
Z \byd \gam\Sha + h E \,,
\eEq
where
$f = \ome (X) \,,\; g = \ome (Y) \,,\; h = \ome (Z) \,.$

 Then, the following equality holds
\bAl\label{expression of d Omega}
	d\Ome(X,Y,Z)
&=
	(i_{E\wed (\Lam\Sha\ten\Lam\Sha)(d\ome)} 
- \tfr12i_{[\Lam,\Lam]}) 		(\alp \wed \bet \wed \gam)
\\ \nonumber
&\quad 
	+ f \, \big(i_{[E,\Lam]} 
+ i_{E \wed (L_E\ome)\Sha}\big)(\bet \wed \gam)
\\ \nonumber
&\quad 
	+ g \, \big(i_{[E,\Lam]} 
+ i_{E \wed (L_E\ome)\Sha}\big)(\gam \wed \alp)
\\ \nonumber
&\quad 
	+ h \, \big(i_{[E,\Lam]} 
+ i_{E \wed (L_E\ome)\Sha}\big)(\alp \wed \bet) \,.
\end{align}
\eLm

\bpf
 Let 
$\wti\alp \,,\, \wti\bet \,,\, \wti\gam$
be the projections of 
$\alp, \bet,\gam$
on
$\sec(\f M, \im (\Ome\Fla)) \sub \sec(\f M, T^*\f M) \,.$

 We have
\beq
	d\Ome(X,Y,Z) =
\eeq
\begin{align*}
&= 
	d\Ome(\alp\Sha + \ome(X) \, E \,, 
	\bet\Sha + \ome(Y) \, E \,, 
	\gam\Sha + \ome(Z) \, E) 
\\
&= 
	d\Ome(\alp\Sha,\bet\Sha,\gam\Sha) 
	+ \ome(X)\, d\Ome(E,\bet\Sha,\gam\Sha) 
	+ \ome(Y)\, d\Ome(\alp\Sha,E,\gam\Sha)
	+ \ome(Z)\, d\Ome(\alp\Sha,\bet\Sha,E) \,.
\end{align*}

 Then, we obtain
\beq
	d\Ome(\alp\Sha, \, \bet\Sha, \, \gam\Sha) =
\eeq
\begin{align*}
&=
	\alp\Sha.\Ome(\bet\Sha,\gam\Sha) 
	+ \bet\Sha.\Ome(\gam\Sha,\alp\Sha)
	+ \gam\Sha.\Ome(\alp\Sha,\bet\Sha)
\\
&
	- \Ome([\alp\Sha,\bet\Sha],\gam\Sha)
	- \Ome([\bet\Sha,\gam\Sha],\alp\Sha)
	- \Ome([\gam\Sha,\alp\Sha],\bet\Sha)
\\[3mm]
&=
	- \alp\Sha.\Lam(\bet,\gam) 
	- \bet\Sha.\Lam(\gam,\alp)
	- \gam\Sha.\Lam(\alp,\bet)
	+ i_{[\alp\Sha,\bet\Sha]}i_{\gam\Sha}\Ome
	+ i_{[\bet\Sha,\gam\Sha]}i_{\alp\Sha}\Ome
	+ i_{[\gam\Sha,\alp\Sha]}i_{\bet\Sha}\Ome
\\[3mm]
&=
	- i_{\alp\Sha}d(\Lam(\bet,\gam)) 
	- i_{\bet\Sha}d(\Lam(\gam,\alp))
	- i_{\gam\Sha}d(\Lam(\alp,\bet))
\\
&
	+ (i_{\alp\Sha}di_{\bet\Sha} - i_{\bet\Sha}di_{\alp\Sha})
		i_{\gam\Sha}\Ome
	+ (i_{\bet\Sha}di_{\gam\Sha} - i_{\gam\Sha}di_{\bet\Sha})
		i_{\alp\Sha}\Ome
	+ (i_{\gam\Sha}di_{\alp\Sha} - i_{\alp\Sha}di_{\gam\Sha})
		i_{\bet\Sha}\Ome
\\
&
	- i_{\alp\Sha \wed \bet\Sha}di_{\gam\Sha}\Ome
	- i_{\bet\Sha \wed \gam\Sha}di_{\alp\Sha}\Ome
	- i_{\gam\Sha \wed \alp\Sha}di_{\bet\Sha}\Ome
\\[3mm]
&=
	\Lam(\alp,d(\Lam(\bet,\gam)))
	+ \Lam(\bet,d(\Lam(\gam,\alp)))
	+ \Lam(\gam,d(\Lam(\alp,\bet)))
\\
&
	- d\ti\alp(\bet\Sha,\gam\Sha)
	- d\ti\bet(\gam\Sha,\alp\Sha)
	- d\ti\gam(\alp\Sha,\bet\Sha)
\\[3mm]
&=
	- i_\Lam d i_\Lam (\alp \wed \bet \wed \gam)
	+ \alp(E)\, d\ome(\bet\Sha,\gam\Sha)
	+ \bet(E)\, d\ome(\gam\Sha,\alp\Sha)
	+ \gam(E)\, d\ome(\alp\Sha,\bet\Sha)
\\[3mm]
&=
	- i_\Lam d i_\Lam (\alp \wed \bet \wed \gam)
	+ \alp(E)\, (\Lam\Sha \ten \Lam\Sha)(d\ome)(\bet,\gam)
\\
&
	+ \bet(E)\, (\Lam\Sha \ten \Lam\Sha)(d\ome)(\gam,\alp)
	+ \gam(E)\, (\Lam\Sha \ten \Lam\Sha)(d\ome)(\alp,\bet)
\\[3mm]
&=
	(i_{E\wed(\Lam\Sha \ten \Lam\Sha)(d\ome)} 
		- i_\Lam d i_\Lam) (\alp \wed \bet \wed \gam)
\\[3mm]
&= 
	(i_{E\wed(\Lam\Sha \ten \Lam\Sha)(d\ome)} 
		- \tfr12 i_{[\Lam,\Lam]}) (\alp \wed \bet \wed \gam)
\,. 
\end{align*}

 Similarly, we obtain
\beq
	d\Ome(\alp\Sha,\bet\Sha,E) =
\eeq
\begin{align*}
&= 
	E.\Ome(\alp\Sha,\bet\Sha) 
	- \Ome([\bet\Sha,E],\alp\Sha)
	- \Ome([E,\alp\Sha],\bet\Sha)
\\
&=
	- E.\Lam(\alp,\bet) 
	+ i_{[\bet\Sha,E]}i_{\alp\Sha} \Ome
	+ i_{[E,\alp\Sha]} i_{\bet\Sha} \Ome
\\
&=
	- E.\Lam(\alp,\bet) 
	+ (i_{\bet\Sha} d i_E - i_E d i_{\bet\Sha} - i_{\bet\Sha \wed E}d)
\ti\alp
	+ (i_E d i_{\alp\Sha}  - i_{\alp\Sha} d i_E - i_{E\wed \alp\Sha}d)
 \ti\bet
\\
&=
	E.\Lam(\alp,\bet) 
	+ i_{\bet\Sha \wed E}d(\alp(E)\,\ome)
	+ i_{E \wed \alp\Sha}d(\bet(E)\,\ome)
\\
&=
	i_Ed i_\Lam (\alp \wed \bet) 
	- \Lam(d(\alp(E)),\bet) 
	+ \Lam(d(\bet(E)),\alp)  
- (\alp(E)\,d\ome(E,\bet\Sha)
	+ (\bet(E)\,d\ome(E,\alp\Sha)
\\
&=
	(i_{[E,\Lam]} 
+ i_{E\wed (L_E\ome)\Sha}) (\alp \wed \bet) \,.
\end{align*}

 Then, the above equalities imply \eqref{expression of d Omega}.
\epf

It is well known \cite{Kir76,Lich78} 
that if
$(\ome, \Ome)$ 
is contact, then
$(E,\Lam)$
is Jacobi.
 Thus, the geometric structures given by dual contact and regular
Jacobi pairs are essentially the same. 
 But we obtain the  equivalence  of structures also for
other types of dual covariant and  contravariant pairs.

\bTh
 The following facts hold:

\smallskip

{\rm (1)}
$(\ome,\Ome)$ 
is an almost--cosymplectic--contact pair 
if and only if 
$(E, \Lam, \ome)$ 
is an almost--coPoisson--Jacobi 3--plet;

\smallskip

{\rm (2)}
$(\ome,\Ome)$ 
is a cosymplectic pair
if and only if 
$(E,\Lam)$ 
is a coPoisson pair;

\smallskip

{\rm (3)}
$(\ome,\Ome)$ 
is a contact pair 
if and only if 
$(E,\Lam)$ 
is a Jacobi pair.
\eTh

\bpf
 Let us consider a point
$x \in \f M \,.$
 All 1--forms on $\f M$ can be obtained, pointwisely, from closed
1--forms. 
 Then, according to the splitting 
\eqref{splitting of X, Y, Z}, 
all vectors
$X, Y, Z \in T_x\f M$
can be obtained, pointwisely, by means of closed forms; conversely,
all closed forms 
$\alp, \bet, \gam \in \sec(\f M, T^*\f M)$
can be obtained, pointwisely, from all vectors above.

 Therefore, from Lemma \ref{Lemma: expression of d Omega}
we deduce the following facts, by means of a pointwise reasoning, by
taking into account the fact that the equality
\eqref{expression of d Omega}
involves the vectors
$X, Y, Z$
and the forms
$\alp, \bet, \gam$
only pointwisely and by considering their arbitrariness at
$x \in \f M \,.$

1) $d\Ome = 0$ 
if and only if 
$[\Lam,\Lam] = 2 E \wed (\Lam\Sha \ten \Lam\Sha)(d\ome)$ 
and $[E,\Lam] = - E\wed
(L_E\ome)\Sha \,,$
i.e. the pair
$(\ome,\Ome)$ 
is almost--cosymplectic--contact if and only if the 3--plet
$(E,\Lam,\ome)$ 
is almost--coPoisson--Jacobi.
 
 2) Moreover, if 
$d\Ome=0$ 
and
$d\ome = 0$ 
then 
$[E,\Lam] = 0$ 
and
$[\Lam,\Lam] = 0 \,,$
i.e.
$(E,\Lam)$ 
is coPoisson. 

 On the other hand, if 
$d\Ome = 0$ 
and $
(E,\Lam)$ 
is coPoisson, then 
$(\Lam\Sha \ten \Lam\Sha)(d\ome) = 0$ 
and 
$(L_E\ome)\Sha = 0 \,,$ 
i.e.
$d\ome (\alp\Sha, \bet\Sha) = 0$ 
and 
$d\ome(E,\alp\Sha) = 0 \,,$ 
for all 1--forms 
$\alp,\, \bet \in \sec(\f M, T^*\f M) \,.$
 Then, from the splitting 
$T\f M = \lang E\rang \drs \im(\Lam\Sha) \,,$
we have
$d\ome = 0$ 
and the pair
$(\ome,\Ome)$ 
is cosymplectic.

 Hence the pair 
$(\ome,\Ome)$ 
is cosymplectic if and only if the pair 
$(E,\Lam)$ 
is coPoisson.

 3) Finally, if 
$d\ome = \Ome \,,$
hence 
$d\Ome = 0 \,,$ 
we have 
$[E,\Lam] = -E \wed (L_E\ome)\Sha = 0$ 
and 
$[\Lam,\Lam] =  2 E \wed (\Lam\Sha \ten \Lam\Sha)(\Ome) =
- 2 E \wed \Lam \,,$
hence the pair 
$(E,\Lam)$ 
is Jacobi.

 On the other hand, if 
$d\Ome = 0$ 
and the pair 
$(E,\Lam)$ 
is Jacobi, then 
$(\Lam\Sha \ten \Lam\Sha)(d\ome) = - \Lam$ 
and
$(L_E\ome)\Sha = 0 \,,$ 
i.e.
$d\ome(\alp\Sha, \bet\Sha) = - \Lam(\alp,\bet)$ 
and 
$d\ome(E, \alp\Sha) = 0 \,,$
hence 
$d\ome = \Ome \,,$
i.e.
the pair 
$(\ome,\Ome)$ 
is contact.

 Thus, the pair 
$(\ome,\Ome)$ 
is contact if and only if the pair 
$(E,\Lam)$ 
is Jacobi.
\epf
\subsection{Darboux's charts}
\label{Darboux's charts}
 First, let us consider an almost--cosymplectic--contact structure 
$(\ome,\Ome) \,.$

\bNt\cite{LibMar87}
 In a neighborhood of each 
$x\in \f M$ 
there exists a local chart (a \emph{Darboux's chart}) 
$(t,x^i,x^{i+n}) \,,$
with
$i = 1, \dots, n \,,$ 
adapted to an almost--cosymplectic--contact structure 
$(\ome,\Ome) \,,$
i.e. such that
\bEq\label{Darboux almost-cosymplectc-contact}
\ome = dt + \sum_{1 \leq i \leq n} 
(\ome^{i} \, dx^i + \ome^{i+n} \, dx^{i+n}) \,,
\qquad
\Ome = \sum_{1 \leq i \leq n}  dx^i \wed dx^{i+n} \,,
\eEq
where 
$\ome^i, \ome^{i+n}\in \map(\f M,\Rn) \,.$

 Indeed, the above almost--cosymplectic--contact pair is cosymplectic
if, for instance, 
$\ome^i = \ome^{i+n}=0$ 
\cite{deLTuy96} 
and contact if, for instance, 
$\ome^i = - x^{i+n}$
and 
$\ome^{i+n} = 0 \,.$
\hfill\ENDE
\eNt

\smallskip

 Then, let us consider an almost--coPoisson--Jacobi structure
$(E,\Lam,\ome) \,.$
 We can find Darboux's charts adapted to this structure,
analogously to the case of almost--cosymplectic--contact structures.

\bLm\label{Lm1a}
 Let 
$\alp, \bet \in \sec(\f M, T^*\f M) \,.$
 Then, we have
\bal
	[E,\alp\Sha] 
&=
	\big(i_E d\alp - \alp(E)\, (L_E\ome)\big)\Sha + 
\Lam(L_E\ome,\alp) \, E \,,
\\
	[\alp\Sha, \bet\Sha] 
&=
	\big(d\Lam(\alp,\bet) 
+ \tfr12 i_{\bet\Sha} d\alp 
+ \tfr12 \alp(E) \, (i_{\bet\Sha}d\ome)
- \tfr12 i_{\alp\Sha} d\bet
- \tfr12 \bet(E) \, (i_{\alp\Sha}d\ome)
\big) \Sha 
\\
&\quad
+ \tfr12 d\ome(\alp\Sha, \bet\Sha) \, E \,.
\end{align*}
\eLm

\bpf
 For each
$h \in \map(\f M, \Rn) \,,$
we have
\bal
	[E,\alp\Sha] . h
&=
	E.(\alp\Sha . h) - \alp\Sha.(E. h) = 
E. \Lam(\alp, dh) - \Lam(\alp, d(E.h))
\\
&=
	i_{[E,\Lam]} (\alp \wed dh) + (i_Ed\alp)\Sha.h 
\\
&=
	- i_{E \wed (L_E\ome)\Sha} (\alp \wed dh) 
		+ (i_E d\alp) \Sha.h
\\
&= 
	- i_E \alp((L_E\ome)\Sha. h) 
+ \Lam(L_E\ome,\alp) (E.h) + (i_Ed\alp) \Sha.h
\\
&=
	(i_Ed\alp - \alp(E) \, (L_E\ome))\Sha.h + \Lam(L_E\ome,\alp) \, E.h
\end{align*}
and
\bal
[\alp\Sha, \bet\Sha].h 
&=
	\alp\Sha.(\bet\Sha.h) - \bet\Sha.(\alp\Sha.h)
\\
&= 
	\Lam(\alp, d\Lam(\bet,dh)) - \Lam(\bet, d\Lam(\alp,dh))
\\
&=
	- \tfr12 i_{[\Lam,\Lam]} (\alp\wed \bet \wed dh) -  
\Lam(dh, d\Lam(\alp,\bet)) +
		\tfr12 (i_{\bet\Sha}d\alp) \Sha.h - 
\tfr12 (i_{\alp\Sha}d\bet) \Sha.h
\\
&=
	(d\Lam(\alp,\bet) + 
\tfr12 i_{\bet\Sha}d\alp 
- \tfr12i_{\alp\Sha}d\bet)\Sha. h + 
		\tfr12 i_{E \wed (\Lam\Sha \ten \Lam\Sha)d\ome} 
(\alp \wed \bet \wed dh)
\\
&=
	(d\Lam(\alp,\bet)
+ \tfr12 i_{\bet\Sha}d\alp 
+ \tfr12 \alp(E)\,(i_{\bet\Sha}d\ome)
- \tfr12i_{\alp\Sha}d\bet
- \tfr12 \bet(E)\,(i_{\alp\Sha}d\ome))\Sha. h 
\\
&\quad + 
		\tfr12 d\ome (\alp\Sha, \bet\Sha) \, E . h \,.
\end{align*}
\vglue-1.6\baselineskip
\epf

\bPr\label{Lm1}
 If
$f,g \in \map(\f M, \Rn) \,,$
then
\bal
	[E,df\Sha] 
&=
- \alp(E) \, (L_E\ome)\Sha
	+ \Lam(L_E\ome,df) \, E \,,
\\
	[df\Sha, dg\Sha] 
&=
	\big(d\Lam(df,dg)
+ \tfr12 E.f \, (i_{dg\Sha}d\ome)
- \tfr12 E.g \, (i_{df\Sha}d\ome)\big)\Sha 
+ \tfr12 d\ome(df\Sha,dg\Sha) \, E \,.
\end{align*}
\ePr

\bpf
 It follows from the above Lemma \ref{Lm1a}, by putting 
$\alp = df$
and
$\bet = dg \,.$
\epf

\bTh\label{Th1a}
 The 
$(2s+1)$-dimensional 
distribution 
$\lang E\rang \drs \im\Lam\Sha$ 
is completely integrable and
$(E,\Lam,\ome)$ 
induces a regular almost--coPoisson--Jacobi structure on the integral
submanifolds of
$\lang E\rang \drs \im\Lam\Sha \,.$
\eTh

\bpf
 By the above Lemma \ref{Lm1a}, the distribution 
$\lang E\rang \drs \im\Lam\Sha$ 
is involutive and of constant rank, so it is completely integrable.

 Let us consider a 
$(2s+1)$-dimensional integral submanifold 
$\iota : \f N \hto \f M$ 
passing through 
$x \in \f M \,.$

 If 
$\ti f, \ti g \in \map(\f N,\Rn) \,,$
then we can extend them
(locally) to 
$f, g \in \map(\f M,\Rn) \,,$ 
such that 
$\ti f = f \com \iota, \, \ti g = g \com \iota \,.$
 Then, we define
$E_{\f N} \in \sec(\f N, T\f N)$ 
and
$\Lam_{\f N} \in \sec(\f N, \Lam^2T\f N)$
by
\beq
E_{\f N}.\ti f = E.f \,,
\quad
\Lam_{\f N}(d\ti f, d\ti g) = 
\Lam(df,dg) = 
(df\Sha).g = 
- (dg)\Sha.f \,.
\eeq
 Indeed, the above
$E_\f N$
and
$\Lam_\f N$ 
depend only on 
$\ti f, \ti g \,,$ 
since they are computed along the integral curves of 
$E, (df)\Sha, (dg)\Sha$ 
through 
$x$ 
and these curves belong to 
$\f N \,.$

 Clearly, $(E_{\f N}, \Lam_{\f N})$ 
satisfy the equalities
\beq
E_{\f N}(\iota^*\alp) = E(\alp) \com \iota \,,
\qquad 
\Lam_{\f N}(\iota^*\alp,\iota^*\bet) = \Lam(\alp,\bet) \com \iota \,,
\qquad
\Al \alp, \bet \in \sec(\f M, T^*\f M) \,.
\eeq

 Then, from the naturality of the Schouten bracket \cite{KolMicSlo93},
we have
\bal
[E_{\f N},\Lam_{\f N}](\iota^*\alp,\iota^*\bet) 
&= 
[E,\Lam](\alp,\bet) \com \iota\,,
\\
[\Lam_{\f N},\Lam_{\f N}](\iota^*\alp,\iota^*\bet,\iota^*\gam) 
	&= 
[\Lam,\Lam](\alp,\bet,\gam) \com \iota\,.
\end{align*}

 Let us set 
$\ome_{\f N} = \iota^*\ome$ 
and
$\Lam\ShaN : T^*\f N \to T\f N: 
\iota^*\alp \mto (\iota^*\alp)\ShaN \byd i_{\iota^*\alp} \Lam_\f N
\,.$

 Then ,
$i_E\ome = 1$
implies 
$i_{E_{\f N}}\ome_{\f N} = 1$ 
and 
$i_\ome\Lam = 0$ 
implies
$i_{\ome_{\f N}}\Lam_{\f N} = 0 \,.$

 Moreover, for each
$\alp, \bet \in \sec(\f M, T^*\f M) \,,$
we have
\beq
\iota^*(L_E\ome) = L_{E_{\f N}}\ome_{\f N} 
\ssep{and}
d\ome_{\f N}\big((\iota^*\alp){\ShaN}, (\iota^*\bet){\ShaN}\big) 
= d\ome(\alp\Sha,\bet\Sha) \com \iota \,.
\eeq

 Then, we have
\beq
[E_{\f N},\Lam_{\f N}]
(\iota^*\alp,\iota^*\bet) =
\eeq
\begin{align*}  
&= 
[E,\Lam](\alp,\bet) \com \iota
\\
&	= 
- (E\wed (L_E\ome)\Sha)(\alp,\bet)\com \iota =
	- E(\alp)\,\Lam(L_E\ome,\bet) \com \iota +
E(\bet)\,\Lam(L_E\ome,\alp) \com \iota 
\\
&=
	- E_{\f N}(\iota^*\alp)\,\Lam_{\f N}(\iota^*(L_E\ome),\iota^*\bet) +
E_{\f N}(\iota^*\bet)\,\Lam_{\f N}(\iota^*(L_E\ome),\iota^*\alp)
\\
&=
	- (E_{\f N} \wed (L_{E_{\f N}}\ome_{\f N}){\ShaN})
(\iota^*\alp,\iota^*\bet) \,.
\end{align*}

 Similarly, we have
\beq
[\Lam_{\f N},\Lam_{\f N}](\iota^*\alp,\iota^*\bet,\iota^*\gam) =
\eeq
\begin{align*}
	&= 
[\Lam,\Lam](\alp,\bet,\gam)\circ \iota
\\
&=
	2 ( E\wed (\Lam\Sha\ten\Lam\Sha)d\ome) (\alp,\bet,\gam )\circ \iota
\\
&=
	2\big(E(\alp)d\ome(\bet\Sha,\gam\Sha) -
E(\bet)d\ome(\alp\Sha,\gam\Sha)
	+ E(\gam)d\ome(\alp\Sha,\bet\Sha)\big)\circ \iota
\\
&=
		2\big(E_{\f N}(\iota^*\alp)d\ome_{\f
N}(\iota^*\bet){\ShaN},(\iota^*\gam){\ShaN}) 
	- E_{\f N}(\iota^*\bet)d\ome_{\f N}((\iota^*\alp){\ShaN},(\iota^*\gam)
{\ShaN})
\\
&\qquad	
	+ E_{\f N}(\iota^*\gam)d\ome_{\f N}((\iota^*\alp){\ShaN},(\iota^*\bet)
{\ShaN})\big)
\\
&=
	2 \, ( E_{\f N} \wed (\Lam_{\f N}\Sha \ten \Lam_{\f N}\Sha)
d\ome_{\f N} )
		(\iota^*\alp,\iota^*\bet,\iota^*\gam ) \,.
\end{align*}	

 Hence, 
$(E_{\f N}, \, \Lam_{\f N} \,, \ome_{\f N})$ 
is a regular almost--coPoisson--Jacobi 3--plet on 
$\f N \,.$
\epf 

\bPr
 In a neighborhood of each 
$x\in \f M$ 
there exists a local chart (a \emph{Darboux's chart}) 
$(W;t,x^i,x^{i+n}) \,,$
with
$i = 1, \dots, n \,,$ 
adapted to the almost--coPoisson--Jacobi 3--plet
$(E, \Lam, \ome)$
i.e. such that
\bAl\label{Darboux almost--coPoisson--Jacobi}
E & = 
\fr{\der}{\der t} \,,
\\ \nonumber
\Lam 
&= 
\sum_{1 \leq i \leq s} 
\fr{\der}{\der x^{i+n}}\wed \fr{\der}{\der x^i}
	- \sum_{1 \leq i \leq s}
\Big(\ome^{i+n}\,\fr{\der}{\der t}\wed \fr{\der}{\der x^i}
	- \ome^{i}\,\fr{\der}{\der t}\wed \fr{\der}{\der x^{i+n}}\Big) \,,
\\ \nonumber
\ome 
&= 
dt + \sum_{1 \leq i \leq n}
\big(\ome^i \, dx^i + \ome^{i+n} \, dx^{i+n}\big) \,,
\end{align}
where 
$\ome^i \,,\; \ome^{i+n} \in \map(\f M,\Rn) \,.$
\ePr

\bpf
 First, let us suppose that
$\Lam$ 
be  of rank 
$2s = 2n$
and let us consider a Darboux's chart adapted to the dual
almost--cosymplectic--contact pair 
$(\ome,\Ome) \,.$
 Then,
from (\ref{Darboux almost-cosymplectc-contact}) we can easily see
that (\ref{Darboux almost--coPoisson--Jacobi}) is satisfied.

 Next, let us suppose that
$2s < 2n \,.$

 Let 
$s = 0 \,.$ 
 Then,
$i_E\ome = 1$ 
implies that there exists a chart 
$(t, x^i, x^{i+n})$
such that
\beq
E = 
\fr{\der}{\der t} \,,
\quad
\Lam  = 0 \,,
\quad
\ome 
= dt + \sum_{1 \leq i \leq n}
\big(\ome^i \, dx^i + \ome^{i+n} \, dx^{i+n}\big) \,.
\eeq

 Let 
$s > 0 \,.$
 Then, let us consider an integral submanifold 
$\f N$ 
of the distribution 
$\lang E\rang \drs \im\Lam\Sha \,.$
 There exists a coordinate neighborhood 
$(W;t, x^i, x^{i+n})$
of each 
$x \in \f N \,,$
with 
$i = 1,\dots,n \,,$ 
such that 
$\f N$ 
is given by 
$x^j = 0, \, x^{j+n}=0 \,,$
with
$j = s+1,\dots, n \,,$ 
and such that the coordinate neighborhood
$(W \cap N;t,x^i,x^{n+i}) \,,$
with
$i = 1,\dots,s \,,$ 
is the Darboux's chart on 
$\f N$ 
adapted to 
$(E_\f N, \Lam_\f N, \ome_\f N) \,.$ 
\epf

\bRm
 It is easy to see that 
$(E,\Lam,\ome)$ 
given by  
(\ref{Darboux almost--coPoisson--Jacobi})
satisfies the conditions for almost--coPoisson--Jacobi 3--plets.
 Indeed, we have
\bEq \label{ELam Darboux}
[E,\Lam] =
\fr{\der}{\der t}\wed \sum_{1 \leq i \leq s} \big(
	- \fr{\der \ome^{i+n}}{\der t} \fr{\der}{\der x^i}
	+ \fr{\der \ome^{i}}{\der t} \fr{\der}{\der x^{i+n}}\big) \,,
\eEq
\bAl\label{LamLam Darboux}
	[\Lam,\Lam] 
&=
	2\fr{\der}{\der t} 
	\wed \Big[\sum^s_{i,j=1} \big(
	\ome^{j+n}\fr{\der \ome^{i+n}}{\der t} 
		+ \fr{\der \ome^{j+n}}{\der x^{i+n}}
	\big) \, \fr{\der}{\der x^i}\wed \fr{\der}{\der x^j}
\\
&\nonumber
	+ \sum^s_{i,j=1}  \big(
		\ome^{i+n}\fr{\der \ome^{j}}{\der t} 
		- \ome^{j}\fr{\der \ome^{i+n}}{\der t} 
		+ \fr{\der \ome^{i+n}}{\der x^{j}}
		- \fr{\der \ome^{j}}{\der x^{i+n}}
	\big) \, \fr{\der}{\der x^i}\wed \fr{\der}{\der x^{j+n}}
\\
&\nonumber
	+ \sum^s_{i,j=1}  \big(
	\ome^{j}\fr{\der \ome^{i}}{\der t} 
		+ \fr{\der \ome^{j}}{\der x^{i}}
	\big) \, \fr{\der}{\der x^{i+n}}\wed \fr{\der}{\der x^{j+n}}
\Big] \,,
\end{align}
\bAl
\Lam\Sha(L_E\ome) 
&= 
\sum_{1 \leq i \leq s} \Big(\fr{\der \ome^{i+n}}{\der t} \,
\fr{\der}{\der x^i}
	- \fr{\der \ome^{i}}{\der t} \, \fr{\der}{\der x^{i+n}} 
	+ (\fr{\der \ome^i}{\der t} \, \ome^{i+n} 
	- \fr{\der \ome^{i+n}}{\der t} \, \ome^i) \, \fr{\der}{\der t}\Big)
\,,
\end{align}
\bEq
	(\Lam\Sha\ten\Lam\Sha)(d\ome) =
\eEq
\bal
&
= 
	\sum_{1 \leq i,j \leq s} ( 
\ome^{i+n}\,\ome^{j+n}\fr{\der \ome^i}{\der t}
		-\ome^{i}\,\ome^{j+n}\fr{\der \ome^{i+n}}{\der t}
		- \ome^{i} \fr{\der \ome^{j+n}}{\der x^{i+n}}
\\
&\qquad\qquad
+ 		\ome^{i+n} \fr{\der \ome^{j+n}}{\der x^{i}}
- \ome^{i+n} \fr{\der \ome^{i}}{\der x^{j+n}}
		+ \ome^{i} \fr{\der \ome^{i+n}}{\der x^{j+n}}) \,
	\fr{\der}{\der t}\wed \fr{\der}{\der x^j}
\\
&
+ \sum_{1 \leq i,j \leq s} ( 
- \ome^{i+n}\,\ome^{j}\fr{\der \ome^i}{\der t}
+ \ome^{i}\,\ome^{j}\fr{\der \ome^{i+n}}{\der t}
+ \ome^{i+n} \fr{\der \ome^{i}}{\der x^{j}}
\\
&\qquad\qquad
		+ \ome^{i+n} \fr{\der \ome^{j}}{\der x^{i}}
- \ome^{i} \fr{\der \ome^{i+n}}{\der x^{j}}
		+ \ome^{i} \fr{\der \ome^{j}}{\der x^{i+n}}) \,
	\fr{\der}{\der t}\wed \fr{\der}{\der x^{j+n}}
	\nonumber
\\
&
	+ \sum_{1 \leq i,j \leq s} ( 
\ome^{j+n} \fr{\der \ome^{i+n}}{\der t}
		+ \fr{\der \ome^{j+n}}{\der x^{i+n}}		) \,
	\fr{\der}{\der x^i}\wed \fr{\der}{\der x^{j}}
\\
&
+ 	\sum_{1 \leq i,j \leq s} ( 
\ome^{i+n} \fr{\der \ome^{j}}{\der t}
		- \ome^{j} \fr{\der \ome^{i+n}}{\der t}
		+ \fr{\der \ome^{i+n}}{\der x^{j}}
		- \fr{\der \ome^{j}}{\der x^{i+n}}) \,
	\fr{\der}{\der x^i} \wed \fr{\der}{\der x^{j+n}}
\\
&
	+ \sum_{1 \leq i,j \leq s} (\ome^{j} \fr{\der \ome^{i}}{\der t}
		+ \fr{\der \ome^{j}}{\der x^{i}}) \,
	\fr{\der}{\der x^{i+n}} \wed \fr{\der}{\der x^{j+n}}
\,.
\end{align*}
\vglue-1.9\baselineskip{ }\hfill\ENDE
\eRm

\smallskip
\bRm
 Let 
$(E,\Lam)$
be a contravariant pair with
$s < n \,.$
 Then,
there exists a 1--form
$\ome$ 
which satisfies 
$i_E\ome = 1$ 
and 
$i_\ome\Lam =0 \,,$
(hence also
$i_\ome(E \wed \Lam^s) = \Lam^s$).
 But such a form is not unique.

 Moreover, the 3--plet
$(E, \Lam, \ome)$
turns out to be almost--coPoisson--Jacobi if and only if the 
equalities 
$[E,\Lam] = - E \wed (L_E\ome)\Sha$ 
and 
$[\Lam,\Lam] = E \wed (\Lam\Sha \ten \Lam\Sha)(d\ome)$
are satisfied.
 We can see it in adapted Darboux's charts; in fact, if the
coordinate expressions of 
$E$ 
and 
$\Lam$ 
are given by \eqref{Darboux almost--coPoisson--Jacobi}, then the
functions 
$\ome^i, \ome^{i+n} \,,$
with
$i = 1, \dots, s \,,$
are given uniquely by 
$\Lam \,,$
but 
$\ome^i, \ome^{i+n} \,,$
with
$i = s+1, \dots, n \,,$
are arbitrary, so
$\ome$ 
is not unique.
\hfill\ENDE
\eRm

\bNt
 The almost--coPoisson--Jacobi 3--plet given in Darboux's charts by
\eqref{Darboux almost--coPoisson--Jacobi}
is coPoisson if, for instance, 
$\ome^i = \ome^{i+n} = 0 \,,$
with
$i = 1,\dots,s \,,$
and is Jacobi if, for instance, 
$\ome^i = - x^{i+n} \,,\; \ome^{i+n} = 0 \,,$
with
$i = 1,\dots,s \,.$
\hfill\ENDE
\eNt
\section{Lie algebra of functions}
\label{Lie algebra of functions}
\setcounter{equation}{0}
 Next, we study the algebras of functions associated with the
geometrical structure given by a pre coPoisson pair. 
\subsection{Poisson algebra of functions}
\label{Poisson algebra of functions}
 First, let us start by considering just a 2--vector
$\Lam \in \sec(\f M, \Lam^2T\f M) \,.$

\bDf
 The \emph{Poisson bracket} of functions 
$f,g \in \map(\f M,\Rn)$ 
is defined as 
\bEq
\{f,g\} \byd i_{df\wed dg}\Lam = i_\Lam (df\wed dg) = \Lam(df,dg) \,.
\eEq
\vglue-1.5\baselineskip{ }\hfill\ENDE
\eDf

\bLm\label{Ji for Pb}
  For each
$f,g,h \in \map(\f M, \Rn) \,,$
we have
\bAl
	\big\{\{f,g\},h\big\}
&
	+ \big\{\{g,h\},f\big\}
	+ \big\{\{h,f\},g\big\} =
 \fr12 i_{[\Lam,\Lam]}(df\wed dg\wed dh)\,.
\end{align}
\eLm

\bpf
 We have
\beq
\big\{\{f,g\},h\big\} = \Lam\big(d(\Lam(df,dg),dh)\big)\,.
\eeq

 Then,
\bal
	i_{[\Lam,\Lam]}(df\wed dg\wed dh)
& =
	2 i_\Lam d i_\Lam (df\wed dg\wed dh)
\\
& = 
	2 i_\Lam d\big(\Lam(df,dg)\,dh - \Lam(df,dh)\, dg + 
\Lam(dg, dh)\, df)
\\
& = 
	2 [\Lam(d\Lam(df,dg),dh) + \Lam(d\Lam(dh,df),dg) + 
\Lam(d\Lam(dg, dh),df)]
\\
& =
	2\Big[\big\{\{f,g\},h\big\}	+ \big\{\{g,h\},f\big\} 	+ 
\big\{\{h,f\},g\big\} 
		\Big]
\,.
\end{align*}
\vglue-1.5\baselineskip
\epf

\bPr
 The following facts are equivalent:
\beq
[\Lam,\Lam] = 0 \,;\ltag{(1}
\eeq
\beq
\big\{\{f,g\}, h\big\}	 + 
\big\{\{g,h\}, f\big\}	 + 
\big\{\{h,f\}, g\big\} = 0 \,,
\qquad
\Al f,g,h \in \map(\f M, \Rn)
\,;\ltag{(2}
\eeq

 (3) the bracket
$\{,\}$ 
is a Lie bracket.
\hfill\ENDE
\ePr

 Thus, a Poisson structure yields a Lie algebra of functions.

\bLm
 The following facts are equivalent:
\beq
[df\Sha,dg\Sha].h = d\{f,g\}\Sha.h \,,
\qquad
\Al f,g,h \in \map(\f M, \Rn) \,;\ltag{(1}
\eeq
\beq
\big\{\{f,g\}, h\big\}	 + 
\big\{\{g,h\}, f\big\}	 + 
\big\{\{h,f\}, g\big\} = 0 \,,
\qquad
\Al f,g,h \in \map(\f M, \Rn)
\,.\ltag{(2}
\eeq
\eLm

\bpf
 We have
\bal
d\{f,g\}\Sha.h 
&=
\Lam\big(\{f,g\}, dh\big) =
\big\{\{f,g\},h\big\} \,,
\\
[df\Sha, dg\Sha]. h 
&= 
df\Sha.dg\Sha. h - dg\Sha.df\Sha. h = 
	\big\{f, \{g,h\}\big\} - \big\{g, \{f,h\}\big\} \,.
\end{align*}

 Then,
\beq
\big(d\{f,g\}\Sha - [df\Sha,dg\Sha]\big).h = 
\big\{\{f,g\},h\big\} +
	\big\{\{g,h\},f\big\} + 
\big\{\{h,f\},g\big\} \,.
\eeq
\vglue-1.5\baselineskip
\epf

\bPr
 The following facts are equivalent:
\beq
[df\Sha,dg\Sha] = d\{f,g\}\Sha \,,
\qquad
\Al f, g \in \map(\f M, \Rn) \,;\ltag{(1}
\eeq

 (2) the map
\beq
\Lam\Sha \com d : \map(\f M,\Rn) \to \sec(\f M,T\f M)
\eeq
is a Lie algebra homomorphism with respect to the Poisson bracket 
of functions and the Lie bracket of vector fields;

 (3) $(\Lam)$
is a Poisson structure, i.e.
$[\Lam,\Lam] = 0 \,.$
\hfill\ENDE
\ePr

\bCr
Let 
$(E,\Lam)$ 
be a coPoisson pair
 then 
$\Lam$ 
defines  a  Poisson algebra of functions.
\hfill\ENDE
\eCr
\subsection{Jacobi algebra of functions}
\label{Jacobi algebra of functions}
 Then, let us consider a contravariant pair
$(E, \Lam) \,.$

\bRm
 If
$(E,\Lam)$ 
is a Jacobi pair with
$s > 0 \,,$
then the Poisson bracket does not satisfy the Jacobi identity.
 In fact, the Jacobi identity turns out to be equivalent to
$E \wed \Lam = 0 \,.$
 But this condition conflicts with our hypothesis
$E \wed \Lam^s \nid 0 \,.$
\hfill\ENDE
\eRm

\bDf
 The \emph{Hamiltonian lift} of a function 
$f\in \map(\f M,\Rn)$
is defined to be the vector field 
\bEq\label{Hamiltonian lift}
X_f \byd i_{df}\Lam - f E = df\Sha - f E\,.
\eEq
\vglue-1.5\baselineskip{ }\hfill\ENDE
\eDf

\bDf
 The \emph{Jacobi bracket} of two functions 
$f,g \in \map(\f M,\Rn)$ 
is defined as 
\bEq\label{Jacobi bracket}
[f,g] \byd \{f,g\} - fE.g + gE.f  = 
\Lam(df,dg) - fE.g + gE.f \,.
\eEq
\vglue-1.5\baselineskip{ }\hfill\ENDE
\eDf

\bLm
 For each
$f,g \in \map(\f M, \Rn) \,,$
we have
\bEq
E.\{f,g\} = \{E.f,g\} + \{f,E.g\} + i_{[E,\Lam]}(df\wed dg)\,.
\eEq
\eLm

\bpf
 We have
\bal
	i_{[E,\Lam]}(df\wed dg)
& = 
	i_E d i_\Lam (df\wed dg) - i_\Lam d i_E(df\wed dg)
\\
& = 
	i_E d i_\Lam (df\wed dg) - i_\Lam d(E.f\,dg - E.g\, df)
\\
& = 
	i_E d i_\Lam (df\wed dg) - i_\Lam (d(E.f)\wed dg) - 
i_\Lam(df\wed d(E.g))
\\
& =
	E.\{f,g\} - \{E.f,g\} - \{f,E.g\}\,. 
\end{align*}
\vglue-1.5\baselineskip
\epf

\bLm\label{L: Jacobi property}
 For each
$f,g,h \in \map(\f M, \Rn) \,,$
we have
\bEq
	\big[[f,g],h\big]
	+ 	\big[[g,h],f\big]
	+ \big[[h,f],g\big] =
\eEq
\beq
=
	(\tfr12 i_{[\Lam,\Lam]} + i_{E\wed \Lam})(df\wed dg\wed dh) 
	+ i_{[E,\Lam]}\big( f\,dg\wed dh  
	+ g \, dh\wed df  
	+ h \, df\wed dg\big)\,. 
\eeq
\eLm

\bpf
 We have
\bal
\big[[f,g],h\big]
&=
	\big\{\{f,g\},h\big\} - \{f,g\}(E.h) - \{f,h\}(E.g) + \{g.h\}(E.f)
\\
&	+ 
h\{E.f,g\} + h\{f, E.g\} - f\{E.g,h\} + g\{E.f, h\} +
	h i_{[E,\Lam]} (df\wed dg)
\\
&+
	f(E.g)(E.h) - g(E.f)(E.h) - hfE^2.g  + hg E^2.f \,.
\end{align*}

 Then,
\beq
	\big[[f,g],h\big] + 
\big[[g,h],f\big] 	+ 
\big[[h,f],g\big]
=
\eeq
\bal
&=
		\big\{\{f,g\},h\big\}
	+ \big\{\{g,h\},f\big\}
	+ \big\{\{h,f\},g\big\} 
	+ \{f,g\} (E.h) + \{g,h\} (E.f) + \{h,f\} (E.g)
\\
&
	+ f \, i_{[E,\Lam]}(dg\wed dh)  
	+ g \, i_{[E,\Lam]}(dh\wed df)  
	+ h \, i_{[E,\Lam]}(df\wed dg)
\\[3mm]
&= 
	\big(\tfr 12 i_{[\Lam,\Lam]}
	+ i_{E\wed \Lam}\big) (df\wed dg\wed dh)
	+ i_{[E,\Lam]}\big( f\,dg\wed dh  
	+ g \, dh\wed df  
	+ h \, df\wed dg\big)\,.
\end{align*}
\vglue-1.5\baselineskip
\epf

\bPr\label{Pr2.11}
\cite{Kir76}
 The Jacobi bracket defines a Lie algebra of functions if and
only if 
$[E, \Lam] = 0 \,$ 
and 
$[\Lam, \Lam] = - 2 E \wed \Lam \,.$

So, a Jacobi pair 
$(E,\Lam)$ 
defines a Lie algebra of functions with respect to the Jacobi
bracket (the \emph{Jacobi algebra} of functions).
\hfill\ENDE\ePr

\bRm
 A coPoisson pair does not define a Lie algebra of functions 
with respect to the Jacobi bracket. 
 Indeed, for a coPoisson pair, we have
\beq
\big[[f,g],h\big] 	+ \big[[g,h],f\big] 	+ \big[[h,f],g\big] =
	i_{E \wed \Lam}(d f\wed dg \wed dh) \,,
\eeq
so, in general, the Jacobi identity is not satisfied.
 Indeed, it is satisfied in the particular case when
$E\wed\Lam = 0 \,,$ 
but this condition conflicts with our hypothesis 
$E \wed \Lam^s \nid 0 \,.$
\hfill\ENDE
\eRm

\bLm\label{L: Hamiltonian lift}
 For each
$f,g,h \in \map(\f M, \Rn) \,,$
we have
\bgt
	([X_f,X_g] - X_{[f,g]}).h =
\\
=
	- (\tfr12 i_{[\Lam,\Lam]} + i_{E \wed \Lam})(df \wed dg \wed dh) 
	-  f \, i_{[E,\Lam]}(dg \wed dh) 
	- g \, i_{[E,\Lam]}(df \wed dh) \,. 
\end{gather*}
\eLm

\bpf
 We have
\beq
	[X_f, X_g].h  =
\eeq
\bal
&= 
	[df\Sha,dg\Sha].h - [df\Sha, gE].h - [fE,dg\Sha].h + [fE,gE].h
\\[3mm]
&=
	d\{f,g\}\Sha.h - \big\{\{f,g\},h\big\} - \big\{\{g,h\},f\big\} 
		-\big\{\{h,f\},g\big\} 
\\
&
	- df\Sha(gE.h) + gE.(df\Sha.h) - fE.(dg\Sha.h)
	+ dg\Sha.(fE.h) 
	+ fE.(gE.h) - g E.(fE.h)
\\[3mm]
&=
	d\{f,g\}\Sha.h - \tfr12 i_{[\Lam,\Lam]}(df \wed dg \wed dh) 
- 2\{f,g\} (E.h) 
\\
&
	- g \{f,E.h\} + gE.\{f,h\} + f \{g,E.h\} - fE.\{g,h\}
	+ f (E.g)(E.h) - g (E.f)(E.h)
\\[3mm]
&=
	d\{f,g\}\Sha.h - \tfr12 i_{[\Lam,\Lam]}(df \wed dg \wed dh)
- 2\{f,g\} (E.h) - f \{E.g,h\} + g \{E.f,h\}
\\
&
- f i_{[E,\Lam]}(dg \wed dh) 
+ g i_{[E,\Lam]}(df \wed dh) 
	+ f (E.g)(E.h) - g (E.f)(E.h)
\end{align*}

 On the other hand,
\bal
	X_{[f,g]}.h
&=
	d\{f,g\}\Sha . h 
- \{f,h\}(E.g) 
+ \{g,h\}(E.f) 
- f\{E.g,h\} 
+ g\{E.f,h\}
\\
&
	- \{f,g\} (E.h) 
+ f(E.g)(E.h) 
- g(E.f)(E.h) \,.
\end{align*}

 Hence,
\bgt
	([X_f,X_g] - X_{[f,g]}).h =
\\
=
	- (\tfr12 i_{[\Lam,\Lam]} + i_{E \wed \Lam})(df \wed dg \wed dh) 
	-  f \, i_{[E,\Lam]}(dg \wed dh) 
	- g \, i_{[E,\Lam]}(df \wed dh) \,. 
\end{gather*}
\vglue-1.5\baselineskip
\epf

\bPr\cite{Lich78}.
 The following facts are equivalent: 
\beq
[X_f,X_g] = X_{[f,g]} \,,
\qquad
\Al f,g \in \map(\f M,\Rn) \,;\ltag{(1}
\eeq

 (2) the Hamiltonian lift of functions with respect to a pair 
$(E,\Lam)$
is a Lie algebra homomorphism with respect to the Jacobi bracket and
the Lie bracket of vector fields;

 (3) the pair
$(E,\Lam)$
is a Jacobi structure, i.e.
$[E,\Lam] = 0$
and
$[\Lam,\Lam] = - 2\, E\wed \Lam \,.$
\ePr

\bpf
 The equivalence follows from Lemma \ref{L: Hamiltonian lift} and the
arbitrariness of the functions
$f,g,h \,.$
\epf
\subsection{Uniqueness of the Jacobi structure}
\label{Uniqueness of the Jacobi structure}
 Now, we  revisit  the well known Proposition \ref{Pr2.11}
\cite{Kir76} in the context of our almost--coPoisson--Jacobi
structures.
 Actually, we prove that an almost--coPoisson--Jacobi 3--plet 
$(E,\Lam,\ome)$ defines a
Lie algebra of functions with respect to the Jacobi bracket if and
only if the pair $(E,\Lam)$ is Jacobi. 

 Let us consider an almost--coPoisson--Jacobi 3--plet 
$(E,\Lam,\ome) \,.$

\bLm\label{L: condition for Jacobi bracket}
 The following facts are equivalent:

(1) for each
$f,g,h \in \map(\f M,\Rn) \,,$
\beq
\big(\tfr 12 i_{[\Lam,\Lam]}
	+ i_{E\wed \Lam}\big) (df \wed dg \wed dh)
	+ i_{[E,\Lam]}\big( f\,dg \wed dh  
	+ g \, dh \wed df  
	+ h \, df \wed dg\big) = 0\,,
\eeq

(2) for each
$f,g \in \map(\f M,\Rn) \,,$
\bEq\label{Eq2.7}
\{f,g\} = - d\ome(X_f,X_g)\,.
\eEq
\eLm

\bpf
  We have
\beq
i_{\Lam\Sha(L_E\ome)} df = - d\ome(E,df\Sha)\,.
\eeq

 Then,
\beq
	\big(\tfr 12 i_{[\Lam,\Lam]}
	+ i_{E \wed \Lam}\big) (df\wed dg \wed dh)
	+ i_{[E,\Lam]}\big(f \,dg \wed dh  
	+ g \, dh \wed df  
	+ h \, df \wed dg\big) =
\eeq
\begin{align*}
&= 
	i_{E \wed 
\big(\Lam + (\Lam\Sha \ten \Lam\Sha)(d\ome)\big)}(df \wed dg \wed dh)
	+ i_{E \wed \Lam\Sha(L_E\ome)} \big(f\,dg \wed dh  
	+ g \, dh \wed df  
	+ h \, df \wed dg\big) 
\\[3mm]
&=
	(E.f)\big(\Lam(dg,dh) + (\Lam\Sha \ten \Lam\Sha)(d\ome)(dg,dh)\big)
\\
&
	+ (E.g)\big(\Lam(dh,df) + (\Lam\Sha \ten \Lam\Sha)(d\ome)(dh,df)\big)
\\
&
	+ (E.h)\big(\Lam(df,dg) + (\Lam\Sha \ten \Lam\Sha)(d\ome)(df,dg)\big)
\\
&
	+ f (E.g) d\ome(E,dh\Sha) - f (E.h) d\ome(E,dg\Sha)
\\
&
	+ g (E.h) d\ome(E,df\Sha) - g (E.f) d\ome(E,dh\Sha)
\\
&
	+ h (E.f) d\ome(E,dg\Sha) - h (E.g) d\ome(E,df\Sha)
\\[3mm]
&= 
	 (E.f)\big(\{g,h\} + d\ome(dg\Sha,dh\Sha) 
		- g d\ome(E,dh\Sha) + h d\ome(E,dg\Sha\big)
\\
&
	+ (E.g)\big(\{h,f\} + d\ome(dh\Sha,df\Sha) 
		- h d\ome(E,df\Sha) + f d\ome(E,dh\Sha\big)
\\
&
	+ (E.h)\big(\{f,g\} + d\ome(df\Sha,dg\Sha) 
		- f d\ome(E,dg\Sha) + g d\ome(E,df\Sha\big)
\\[3mm]
&=
	 (E.f)\big(\{g,h\} + d\ome(dg\Sha - g E,dh\Sha - h E) \big)
\\
&
	+ (E.g)\big(\{h,f\} + d\ome(dh\Sha - h E,df\Sha - f E) \big)
\\
&
	+ (E.h)\big(\{f,g\} + d\ome(df\Sha - f E,dg\Sha - g E) \big) \,. 
\end{align*}
\vglue-1.5\baselineskip
\epf

\bPr
 The almost--coPoisson--Jacobi structure 
$(E, \, \Lam, \, \ome)$
yields a Lie algebra of functions with respect to the Jacobi bracket
if and only if the Poisson bracket satisfies
(\ref{Eq2.7}).
\ePr

\bpf
 It follows from the above Lemma \ref{L: condition for Jacobi bracket}
end from Lemma \ref{L: Jacobi property}.
\epf

\bCr
 A Jacobi pair
$(E,\Lam)$
yields a Lie algebra with respect to the Jacobi bracket.

 A coPoisson pair 
$(E,\Lam)$ 
yields a Lie algebra with respect to the Jacobi bracket if and only
if
$\Lam = 0 \,.$
\eCr

\bpf
 Let
$(E,\Lam)$
be a Jacobi pair.
 Then, for each
$\alp, \bet \in \sec(\f M, T^*\f M) \,,$
we have
\beq
d\ome(\alp\Sha, \bet\Sha) = - \Lam(\alp,\bet)
\ssep{and}
d\ome(E,\alp\Sha) = 0 \,,
\eeq
hence, for each
$f,g \in \map(\f M, \Rn) \,,$
we obtain
\bal
\{f,g\} 
&\byd 
\Lam(df,dg) =
-d\ome(df\Sha,dg\Sha) =
-d\ome(df\Sha,dg\Sha) + f\, d\ome(E,dg\Sha) - g\, d\ome(df\Sha,E) 
\\
& 
\byd - d\ome(X_f,X_g) \,,
\end{align*}
hence condition \eqref{Eq2.7} is satisfied.

 Let
$(E,\Lam)$
be a coPoisson pair.
 Then, we have
$d\ome = 0 \,,$
hence condition \eqref{Eq2.7} is satisfied if and only if
$\{f,g\} = 0 \,,$
i.e. if and only if 
$\Lam = 0 \,.$
\epf

\bTh
 An almost--coPoisson--Jacobi 3--plet 
$(E,\Lam,\ome)$ 
yields a Lie algebra of functions with respect to the Jacobi
bracket if and only if the pair $(E,\Lam)$ is Jacobi. 
\eTh

\bpf
 It is sufficient to prove that \eqref{Eq2.7} implies that the
pair
$(E,\Lam)$ 
is Jacobi.

 We can prove it in a local chart.

 In a Darboux's cchart adapted to an almost--coPoisson--Jacobi
3--plet 
$(E,\Lam,\ome)$ 
according to \eqref{Darboux almost--coPoisson--Jacobi} we have
\bMl
	X_f = 
	\Big(- f + \sum_{1 \leq i \leq s}
\big( \ome^{i+n}\fr{\der f}{\der x^i}
		- \ome^i \fr{\der f}{\der x^{i+n}}
	\big)\Big)\fr{\der}{\der t}
\\
	+ \sum_{1 \leq i \leq s} \big(
		\fr{\der f}{\der x^{i+n}} 
		- \ome^{i+n} \fr{\der f}{\der t} 
	\big)\fr{\der}{\der x^{i}}
	+ \sum_{1 \leq i \leq s} \big(
		\ome^i \fr{\der f}{\der t} 
		- \fr{\der f}{\der x^i} 
	\big)\fr{\der}{\der x^{i+n}} \,.
\end{multline}

 Then,
\beq
	d\ome(X_f,X_g) =
\eeq
\begin{align*}
&=
		\big(f\fr{\der g}{\der t} - g \fr{\der f}{\der t}\big) \,.
	\sum_{i=1}^s \big(\fr{\der \ome^i}{\der t}\ome^{i+n} 
		- \fr{\der \ome^{i+n}}{\der t}\ome^i\big)
\\
&+
	\sum_{1 \leq i \leq s} \big(f\fr{\der g}{\der x^i} 
- g\fr{\der f}{\der x^i}\big) \fr{\der\ome^{i+n}}{\der t}
	- \sum_{i=1}^s \big(f\fr{\der g}{\der x^{i+n}} 
- g\fr{\der f}{\der x^{i+n}}\big)
		\fr{\der\ome^{i}}{\der t}
\\
&+
	\sum_{1 \leq i,j \leq s} \big( 
		\fr{\der f}{\der t}\fr{\der g}{\der x^i}
		-\fr{\der g}{\der t}\fr{\der f}{\der x^i}
	\big) \,.
	\Big(
		\ome^{j+n}\ome^{i+n}\fr{\der \ome^j}{\der t}
		- \ome^{j}\ome^{i+n}\fr{\der \ome^{j+n}}{\der t}
\\
&\qquad\qquad
		+ \ome^{j+n}\fr{\der \ome^{i+n}}{\der x^j}
		- \ome^{j+n}\fr{\der \ome^{j}}{\der x^{i+n}}
		+ \ome^{j}\fr{\der \ome^{j+n}}{\der x^{i+n}}
		- \ome^{j}\fr{\der \ome^{i+n}}{\der x^{j+n}}
	\Big)
\\
&+
	\sum_{1 \leq i,j \leq s} \big( 
		\fr{\der f}{\der t}\fr{\der g}{\der x^{i+n}}
		-\fr{\der g}{\der t}\fr{\der f}{\der x^{i+n}}
	\big) \,.
	\Big(
		\ome^{j}\ome^{i}\fr{\der \ome^{j+n}}{\der t}
		- \ome^{j+n}\ome^{i}\fr{\der \ome^{j}}{\der t}
\\
&\qquad\qquad
		+ \ome^{j+n}\fr{\der \ome^{j}}{\der x^i}
		- \ome^{j+n}\fr{\der \ome^{i}}{\der x^{j}}
		+ \ome^{j}\fr{\der \ome^{j+n}}{\der x^{i+n}}
		- \ome^{j}\fr{\der \ome^{j+n}}{\der x^{i}}
		+ \ome^{j}\fr{\der \ome^{i}}{\der x^{j+n}}
		- \ome^{j+n}\fr{\der \ome^{j}}{\der x^{i+n}}
	\Big)
\\
&+
\sum_{1 \leq i,j \leq s} 
\fr{\der f}{\der x^i}\fr{\der g}{\der x^j}\big(
	\ome^{j+n}\fr{\der\ome^{i+n}}{\der t} 
	- \ome^{i+n}\fr{\der\ome^{j+n}}{\der t}
	+ \fr{\der\ome^{i+n}}{\der x^{j+n}}  
	- \fr{\der\ome^{j+n}}{\der x^{i+n}} 
	\big)
\\
&+
	\sum_{1 \leq i,j \leq s} \big( 
		\fr{\der f}{\der x^i}\fr{\der g}{\der x^{j+n}}
		-\fr{\der g}{\der x^i}\fr{\der f}{\der x^{j+n}}
	\big) \,.
	\Big(
		\ome^{i+n}\fr{\der \ome^{j}}{\der t}
		- \ome^{j}\fr{\der \ome^{i+n}}{\der t}
		+ \fr{\der \ome^{i+n}}{\der x^j}
		- \fr{\der \ome^{j}}{\der x^{i+n}}
	\Big)
\\
&+
	\sum_{1 \leq i,j \leq s}^s  
		\fr{\der f}{\der x^{i+n}}\fr{\der g}{\der x^{j+n}}
	\Big(
		\ome^{j}\fr{\der \ome^{i}}{\der t}
		- \ome^{i}\fr{\der \ome^{j}}{\der t}
		+ \fr{\der \ome^{j}}{\der x^i}
		- \fr{\der \ome^{i}}{\der x^{j}}
	\Big) \,.
\end{align*}

 On the other hand,
\bml
	\{f,g\} =
	\sum_{1 \leq i \leq s} \Big(
		\fr{\der f}{\der x^{i+n}}\fr{\der g}{\der x^{i}}
		- \fr{\der g}{\der x^{i+n}}\fr{\der f}{\der x^{i}}
\\
		- \ome^{i+n} \big(
			\fr{\der f}{\der t}\fr{\der g}{\der x^{i}}
			- \fr{\der g}{\der t}\fr{\der f}{\der x^{i}}
		\big)
		+ \ome^{i} \big(
			\fr{\der f}{\der t}\fr{\der g}{\der x^{i+n}}
			- \fr{\der g}{\der t}\fr{\der f}{\der x^{i+n}}
		\big)\Big)
\end{multline*}

 Now, if we assume 
$\{f,g\} = - d\ome(X_f,X_g) \,,$ 
then we obtain the following system of partial differential
equations, by comparing the above expressions, for all 
$i,j = 1,\dots, s \,,$
\bgt
	0 =
	\fr{\der \ome^{i+n}}{\der t} \,,
\qquad
	0 = \fr{\der \ome^{i}}{\der t} \,,
\qquad
	0 = 	\sum_{1 \leq j \leq s} 
	\Big(
		 \ome^{j}\fr{\der \ome^{j+n}}{\der x^{i+n}}
		- \ome^{j+n}\fr{\der \ome^{j}}{\der x^{i+n}}
	\Big)\,,
\\
	0 =
	\Big(
		 \fr{\der \ome^{j}}{\der x^i}
		- \fr{\der \ome^{i}}{\der x^{j}}
	\Big) \,,
\qquad
	0 =
	\Big(
	 \fr{\der\ome^{i+n}}{\der x^{j+n}}  
	- \fr{\der\ome^{j+n}}{\der x^{i+n}} 
	\Big) \,,
\qquad
	\del^i_j =	  
	\Big(
		 \fr{\der \ome^{i+n}}{\der x^j}
		- \fr{\der \ome^{j}}{\der x^{i+n}}
	\Big)\,.
\end{gather*}

 Now, if we use the above identities, then we obtain 
\bal
	[E,\Lam]
&=
	0
\\
	[\Lam,\Lam] 
&= 
	2\fr{\der}{\der t} 
	\wed \Big(\sum_{1 \leq i,j \leq s} 
	\fr{\der \ome^{j+n}}{\der x^{i+n}}
	 \, \fr{\der}{\der x^i}\wed \fr{\der}{\der x^j}
	+ \sum_{1 \leq i \leq s}  \fr{\der}{\der x^i}\wed \fr{\der}{\der
x^{i+n}}
\\
&\qquad\qquad
	+ \sum_{1 \leq i,j \leq s}  \fr{\der \ome^{j}}{\der x^{i}}
	\, \fr{\der}{\der x^{i+n}}\wed \fr{\der}{\der x^{j+n}}
\Big)
\\
&=
	-2 E \wed \Lam \,.
\end{align*}

So, 
$(E,\Lam)$ 
is a Jacobi pair.
\epf
\section{Examples: dynamical structures}
\label{Examples}
 As examples of the geometric structures analysed above, now we
discuss the dynamical structures arising on the phase spase of
a spacetime in classical relativistic theories.
 We consider the relativistic Galilei and the Einstein spacetimes,
emphasizing the analogies and the differences between the two cases.

\smallskip

 In order to make our theory explicitly independent from units of
measurement, we introduce the ``spaces of scales'' \cite{JanModVit07}.
 Roughly speaking, a space of scales
$\B S$
has the algebraic structure of
$\Rn^+$
but has no distinguished ``basis".
 We can naturally define the tensor product of spaces of scales and
the tensor product of spaces of scales and vector spaces.
 We can also naturally define rational tensor powers
$\B S^{m/n}$
of a space of scales
$\B S \,.$
 Moreover, we can make a natural identification
$\B S^* \seq \B S^{-1} \,.$

 The basic objects of our theory (the metric field, the phase 2--form,
the phase 2--vector, etc.) will be valued into \emph{scaled} vector
bundles, that is into vector bundles multiplied tensorially with
spaces of scales.
 In this way, each tensor field carries explicit information on its
``scale dimension".
 Actually, we assume the following basic spaces of scales:
the space of \emph{time intervals}
$\B T \,,$  
the space of \emph{lengths}
$\B L $ 
and the space of \emph{masses} 
$\B M \,.$
 Moreover, we consider the following ``universal scales": 
the \emph{speed of light}
$c \in \B T^{-1} \ten \B L$ 
and the \emph{Planck constant} 
$\h \in \B T^* \ten \B L^2 \ten \B M \,.$

 A \emph{time unit}
is defined to be an element
$u_0 \in \B T \,,$
or, equivalently, its dual
$u^0 \in \B T^* \,,$ 
\subsection{Galilei spacetime}
\label{Galilei spacetime}
 First, we study the geometrical structures arising on the phase
space of a Galilei spacetime
\cite{JadJanMod98, JanMod99, JanMod02, Tra66}.
\subsubsection{Spacetime}
\label{G: Spacetime}
 We assume \emph{absolute time} to be an affine
1--dimensional space
$\f T$
associated with the vector space
$\baB T \byd \B T \ten \Rn \,.$

 We assume spacetime to be an oriented $(3+1)$-dimensional fibred
manifold
$\f E$
equipped with a \emph{time fibring}
$t : \f E \to \f T \,.$

 A \emph{spacetime chart} is defined to be a chart
$(x^\lam) \eqv (x^0, x^i)$
of
$\f E \,,$
adapted to the orientation, to the fibring, to the affine structure of
$\f T$
and to a time unit
$u_0 \,.$
Greek indices will span all spacetime coordinates and Latin indices
will span the fibre coordinates.
 In the following, we shall always refer to spacetime charts.
 The induced local bases of
$T\f E$
and
$T^*\f E$
are denoted, respectively, by
$(\der_\lam)$
and
$(d^\lam) \,.$

 The vertical restriction of forms will be denoted by the
``check" symbol
$^\wch{\,} \,.$

 The differential of the time fibring is the scaled 1--form
$dt : \f E \to \baB T \ten T^*\f E \,,$
with coordinate expression
$dt = u_0 \ten d^0 \,.$

\smallskip

 We assume spacetime to be equipped with a \emph{scaled spacelike
Riemannian metric}
$g : \f E \to \B L^2 \ten (V^*\f E \ten V^*\f E) \,.$
 The contravariant metric is denoted by
$\ba g : \f E \to \B L^{-2} \ten (V\f E \ten V\f E) \,.$

 We have the coordinate expressions
\bat{3}
g
&=
g_{ij} \, \ch d^i \ten \ch d^j \,,
&&\ssep{with}
g_{ij}
&&\in \map(\f E, \, \B L^2 \ten \Rn) \,,
\\
\ba g
&=
g^{ij} \, \der_i \ten \der_j \,,
&&\ssep{with}
g^{ij}
&&\in \map(\f E, \, \B L^{-2} \ten \Rn) \,.
\end{alignat*}
\subsubsection{Phase space}
\label{G: Phase space}
 A \emph{motion} is defined to be a section
$s :\f T \to \f E \,.$
 The \emph{1st differential} of a motion
$s$
is defined to be the map
$ds : \f T \to \B T^* \ten T\f E \,.$
 We have
$dt (ds) = 1 \,.$

 We assume as \emph{phase space} the 1st jet space
$J_1\f E$
of motions.

 A space time chart
$(x^\lam)$
induces naturally a chart
$(x^\lam, x^i_0)$
on
$J_1\f E \,.$

 The \emph{velocity} of a motion
$s$
is defined to be its 1st jet
$j_1s : \f T \to J_1\f E \,.$

 We define the \emph{contact map} to be the unique fibred morphism
$\K d : J_1\f E \to \B T^* \ten T\f E$
over
$\f E$
such that
$\K d \com j_1s = ds \,,$
for each motion
$s \,.$
 We have
$\K d \con dt = 1 \,.$
 The coordinate expression of
$\K d$
is
\beq
\K d = u^0 \ten \K d_0 \eqv u^0 \ten (\der_0 + x^i_0 \, \der_i) \,.
\eeq

 The map
$\K d$
is injective.
 Accordingly, the 1st jet space can be naturally identified with the
subbundle
$J_1\f E \sub \B T^* \ten T\f E \,,$
of scaled vectors which project on
$\f 1 : \f T \to \baB T^* \ten \baB T \,.$
 Thus, the bundle
$J_1\f E \to \f E$
turns out to be affine and associated with the vector bundle
$\B T^* \ten V\f E \,.$
 Indeed, 
$J_1\f E \sub \B T^*\ten T\f E$
is the fibred submanifold over
$\f E$
characterised by the constraint
$\dt x^0_0 = 1 \,.$

 We define also the \emph{complementary contact map}
$\tht \byd 1 - \K d \com dt : J_1\f E \to T^*\f E \ten V\f E \,.$
 The coordinate expression of
$\tht$
is
\beq
\tht = \tht^i \ten \der_i \eqv (d^i - x^i_0 \, d^0) \ten \der_i \,.
\eeq
\subsubsection{Vertical bundle of the phase space}
\label{G: Vertical bundle of the phase space}
 Let
$V_0J_1\f E \sub VJ_1\f E \sub TJ_1\f E$
be the vertical tangent subbundle over
$\f E$
and the vertical tangent subbundle over
$\f T \,,$
respectively.
 The affine structure of the phase space yields the equality
$V_0J_1\f E = J_1\f E \ucar{\f E} (\B T^* \ten V\f E) \,,$
hence the natural map
$\nu : J_1\f E \to \B T \ten (V^*\f E \ten V_0J_1\f E) \,,$
with coordinate expression
$\nu = u_0 \ten \ch d^i \ten \der^0_i \,.$
\subsubsection{Spacetime connections}
\label{G: Spacetime connections}
 We define a \emph{spacetime connection} to be a torsion free linear
connection 
$K : T\f E \to T^*\f E \ten TT\f E$
of the bundle
$T\f E \to \f E \,.$
 Its coordinate expression is of the type
\beq
K =
d^\lam \ten
(\der_\lam + K\col\lam\mu\nu \, \dot x^\nu \, \dot \der_\mu)\,,
\ssep{with} 
K\col\lam\mu\nu = K\col\nu\mu\lam \in \map(\f E,\Rn) \,.
\eeq

 A spacetime connection 
$K$ 
is said to be  \emph{time preserving} if it preserves the time
fibring, i.e. if
$\nab dt = 0 \,.$
 In coordinates, this reads
$K\col\lam 0 \mu = 0 \,.$

 A time preserving spacetime connection 
$K$ 
is said to be  \emph{metric} if it preserves the metric 
$g \,,$ 
i.e. if
$\nab g = 0 \,.$ 
 In coordinates, it reads
\bal
K\col 0i0
&= - g^{ij} \, 2 \, \phi_{0,0j} \,,
\\
K\col 0ih = K\col hi0
&= - \tfr12 g^{ij} \, 
(2 \, \phi_{0,hj} + \der_0 g_{hj}) \,,
\\
K\col kih
= K\col hik
&= - \tfr12 g^{ij} \,
(\der_h g_{jk} + \der_k g_{jh} - \der_j g_{hk}) \,,
\end{align*}
where
$\phi \in \sec(\f E, \B T^* \ten \B L^2 \ten \Lam^2T^*\f E)$
is a scaled spacetime 2--form (which depends on 
$K$ 
and on the chosen chart).

 The vertical restriction of of a metric spacetime connection 
$K$ 
is just the Levi Civita connection of the spacetime fibres. 

 A spacetime connection 
$K$ 
is said to be a \emph{Galilei connection} if it is time preserving,
metric and such that its curvature tensor
$R$
fulfills a symmetry condition which in coordinates reads
$R_{\lam}{}^i{}_\mu{}^j = R_{\mu}{}^j{}_\lam{}^i \,,$
where 
$R_{\lam}{}^i{}_\mu{}^j \byd g^{jp} \, R\col\lam{i}{\mu p} \,.$
\subsubsection{Phase connections}
\label{G: Phase connections}
 We define a \emph{phase connection} to be a connection of the
bundle
$J_1\f E \to \f E \,.$

 A phase connection can be represented, equivalently, by a tangent
valued form
$\Gam : J_1\f E \to T^*\f E \ten T J_1\f E \,,$
which is projectable over
$\1 : \f E \to T^*\f E \ten T\f E \,,$
or by the complementary vertical valued form
$\nu[\Gam] :
J_1\f E \to T^*J_1\f E \ten VJ_1\f E \,,$
respectively, with coordinate expressions
\beq
\Gam =
d^\lam \ten(\der_\lam + \Gam\Ga\lam i \, \der_i^0) \,,
\quad
\nu[\Gam] =
(d^i_0 - \Gam\Ga\lam i \, d^\lam) \ten \der^0_i \,,
\ssep{with}
\Gam\Ga\lam i \in \map(J_1\f E, \, \Rn) \,.
\eeq

 The coordinate expression of an affine phase connection
$\Gam$
is
$\Gam\Ga \lam i = 
\Gam\Gaa \lam i p \, x^p_0 + \Gam\Gaa \lam i 0 \,.$

 We can prove \cite{JanMod96} that there is a natural bijective map
$\chi : K \mto \Gam$
between time preserving linear spacetime connections
$K$
and affine phase connections
$\Gam \,,$
with coor\-dinate expression
$\Gam\Gaa \lam i \mu = K\col \lam i \mu \,.$
\subsubsection{Dynamical phase connection}
\label{G: Dynamical phase connection}
 The space of 2--jets of motions
$J_2\f E$
can be naturally regarded as the affine subbundle
$J_2\f E \sub \B T^* \ten TJ_1\f E \,,$
which projects on
$\K d : J_1\f E \to \B T^* \ten T\f E \,.$

 A \emph{dynamical phase connection} is defined to be a
2nd--order connection, i.e. a section
$\gam : J_1\f E \to J_2\f E \,,$
or, equivalently, a section
$\gam : J_1\f E \to \B T^* \ten TJ_1\f E \,,$
which projects on
$\K d \,.$

 The coordinate expression of a dynamical phase connection is of the
type
\beq
\gam = c \, \alp^0 \,
(\der_0 + x^i_0 \, \der_i + \gam\ga i \, \der^0_i) \,,
\ssep{with}
\gam\ga i \in \map(J_1\f E, \, \Rn) \,.
\eeq

 If
$\gam$
is a dynamical phase connection, then we have
$\gam \con dt = 1 \,.$

 The contact map
$\K d$
and a phase connection
$\Gam$
yield the section
$\gam \eqv \gam[\K d, \Gam] \byd 
\K d\con\Gam : J_1\f E \to \B T^* \ten TJ_1\f E \,,$
which turns out to be a dynamical phase connection, with coordinate
expression
\beq
\gam\ga i = \Gam\Ga 0i + \Gam\Ga ji \, x^j_0 \,.
\eeq

 In particular, a time preserving spacetime connection
$K$
yields the dynamical phase connection
$\gam \byd \gam[\K d, K] \byd \K d \con \chi(K) \,,$
with coordinate expression
\beq
\gam^i_{00} = 
K\col hik \, x^h_0 \, x^k_0 + 2 \, K\col hi0 \, x^h_0 + K\col 0i0 \,.
\eeq
\subsubsection{Phase 2--form and 2--vector}
\label{G: Phase 2--form and 2--vector}
 The metric
$g$
and a phase connection
$\Gam$
yield the scaled 2--form 
$\Ome \,,$ 
called \emph{(scaled) phase
2--form}, and  the scaled vertical 2--vector 
$\Lam \,,$ 
called \emph{(scaled) phase 2--vector},
\bal
\Ome 
&= 
\Ome[g,\Gam] \byd  g \con \big(\nu[\Gam] \wed \tht \big) :
J_1\f E \to \B T^*\ten \B L^2\ten \Lam^2T^*J_1\f E \,,
\\
\Lam 
&=
\Lam[g,\Gam] \byd \ba g \con (\Gam \wed \nu)
: J_1\f E \to \B T \ten \B L^{-2}\ten \Lam^2VJ_1\f E \,,
\end{align*}
with coordinate expressions
\bal
\Ome[g,\Gam] 
&=
g_{ij} \, u^0 \ten (d^i_0 - \Gam\Ga \lam i \, d^\lam) \wed
(d^j - x^j_0 \, d^0) \,,
\\
\Lam[g,\Gam] 
&=
g^{ij} \, u_0 \ten
\big(\der_i + \Gam\Ga ih \, \der^0_h\big) \wed \der^0_j \,.
\end{align*}

 We can easily see that 
$dt \wed \Ome^3 \not\equiv 0\,$
and  
$\gam \wed \Lam^3 \not\equiv 0\,.$

 There is a unique dynamical phase connection
$\gam \,,$
such that
$\gam \con \Ome[g, \Gam] = 0 \,.$
Namely,
$\gam = \gam[\K d, \Gam] \,.$

 In particular, a metric spacetime connection
$K$
yields the (scaled) phase 2--form 
$\Ome \eqv \Ome[g,K]\byd \Ome[g,\chi(K)]$ 
and the (scaled) phase 2--vector
$\Lam \eqv \Lam[g,K]\byd \Lam[g,\chi(K)]$
with coordinate expressions
\bal
\Ome 
&= - 
g_{ij} \, (d^i - x^i_0 \, d^0) \wed d^j_0
+ \big(\tfr12 \der_j g_{hk}\, x^h_0 \, x^k_0
+ \der_0 g_{hj} \, x^h_0 + \phi_{0, 0j} \big) \, d^0 \wed d^j
\\
&\qquad\qquad\qquad\qquad\quad\quad\;\,
+ \big(\tfr12 (\der_i g_{hj} - \der_j g_{hi}) \, x^h_0
+ \tfr12 \phi_{0,ij} \big) \, d^i \wed d^j \,,
\\
\Lam 
&= 
g^{ij} \, \der_i \wed \der^0_j
- \tfr12 g^{ih} \, g^{jk}
\, \big((\der_k g_{lr} - \der_h g_{lk}) \, x^l_0
+ \phi_{0,kh} \big) \, \der_i^0 \wed \der^0_j \,,
\end{align*}
\subsubsection{Dynamical structures of the phase space}
\label{G: Dynamical structures of the phase space}
 We have the following result \cite{JadJanMod98,JanMod99}.

\bTh
 Let us consider a spacetime connection
$K$
and the induced objects
$\Gam \byd \chi(K) \,,$
$\gam \byd \gam[\K d, \Gam] \,,$
$\Ome \byd \Ome[g, \Gam]$
and
$\Lam \byd \Lam[g, \Gam] \,.$
 Then, the following assertions are equivalent.

\smallskip

{\rm (1)} 
$K$ 
is a Galilei connection.

\smallskip

{\rm (2)} 
$\Ome$ 
is closed, i.e. 
$(-dt, \Ome)$ 
is a scaled cosymplectic pair.

\smallskip

{\rm (3)} 
$[\gam, \Lam] = 0$ 
and 
$[\Lam, \Lam] = 0 \,,$ 
i.e.
$(- \gam, \Lam)$ 
is a scaled (regular) coPoisson pair.

\smallskip

 Moreover, the cosymplectic pair 
$(-dt, \Ome)$ 
and the coPoisson pair
$(-\gam, \Lam)$ 
are mutually dual.
\hfill\ENDE
\eTh

\bRm
 If
$K$
is a time preserving spacetime connection, then the induced pairs 
$(-dt,\Ome[g,K])$ 
and
$(-\gam[\K d, K], \, \Lam[g,K])$ 
are scaled.

 On the other hand, some results of the general theory of geometrical
structures developed in the first two sections requires unscaled
pairs.

 Indeed, if we refer to a particle of mass
$m \in \B M \, $ 
and consider the universal scales
$\h \in \B T^{-1} \ten \B L^2 \ten \B M$
and
$c \in \B T^{-1} \ten \B L \,,$
then we obtain unscaled pairs in the following natural way.

 We have the unscaled spacetime 1--form
\beq
\tfr{m\,c^2}{\h} \, dt : \f E \to T^*\f E \,.
\eeq

 Moreover, the rescaled contact map
$\K D \byd \tfr{\h}{m\, c^2} \, \K d : J_1\f E \to T\f E$
yields the unscaled phase vector field
\beq
\gam \eqv \gam[\K D,K] = 
\tfr{\h}{m\, c^2} \, \gam[\K d,K] : \f E \to TJ_1\f E \,.
\eeq

 Furthermore, the rescaled metric
$G \byd \tfr m{\h} \, g : \f E \to \B T \ten V^*\f E \ten V^*\f E$
yields the unscaled phase 2--form and phase 2--vector
\bal
\Ome \eqv \Ome[G,K] 
&= 
\tfr m{\h} \, \Ome[g,K] : J_1\f E \to \Lam^2T^*J_1\f E \,,
\\
\Lam \eqv \Lam[G,K] 
&= 
\tfr \h{m} \, \Lam[g,K] : J_1\f E \to \Lam^2TJ_1\f E \,.
\end{align*}

 Thus, if
$K$
is a Galilei spacetime connection, then
$(-\tfr{m\,c^2}{\h} \, dt, \, \Ome)$
and 
$(-\tfr{\h}{m\, c^2} \, \gam, \, \Lam)$ 
turn out to be mutually dual unscaled cosymplectic and coPoisson
pairs of the phase space.

 Indeed, the Plank constant does not play any direct role in
classical mechanics; nevertheless, such a scale is necessary for
getting unscaled objects as above.
\hfill\ENDE
\eRm
\subsection{Einstein spacetime}
\label{Einstein spacetime}
\setcounter{equation}{0}
 Then, we study the geometrical structures arising on the phase
space of an Einstein spacetime
\cite{JanMod95,JanMod06}.
\subsubsection{Spacetime}
\label{E: Spacetime}
 We assume \emph{spacetime} to be an oriented 4--dimensional
manifold
$\f E$
equipped with a scaled Lorentzian metric
$g : \f E \to \B L^2 \ten (T^*\f E \ten T^*\f E) \,,$
with signature
$(-+++) \,;$
we suppose spacetime to be time oriented.
 The contravariant metric is denoted by
$\ba g : \f E \to \B L^{-2} \ten (T\f E \ten T\f E) \,.$

 A \emph{spacetime chart} is defined to be a chart
$(x^\lam) \eqv (x^0, x^i) \in \map(\f E, \, \Rn \car \Rn^3)$
of
$\f E \,,$
which fits the orientation of spacetime and such that the vector field
$\der_0$
is timelike and time oriented and the vector fields
$\der_1, \der_2, \der_3$
are spacelike.
 Greek indices
$\lam, \mu, \dots$
will span spacetime coordinates, while Latin indices
$i, j, \dots$
will span spacelike coordinates.
 In the following, we shall always refer to spacetime charts.
 The induced local bases of
$T\f E$
and
$T^*\f E$
are denoted, respectively, by
$(\der_\lam)$
and
$(d^\lam) \,.$
 We have the coordinate expressions
\bat{3}
g
&=
g_{\lam\mu} \, d^\lam \ten d^\mu \,,
&&\ssep{with}
g_{\lam\mu}
&&\in \map(\f E, \, \B L^2 \ten \Rn) \,,
\\
\ba g
&=
g^{\lam\mu} \, \der_\lam \ten \der_\mu\,,
&&\ssep{with}
g^{\lam\mu}
&&\in \map(\f E, \, \B L^{-2} \ten \Rn) \,.
\end{alignat*}
\subsubsection{Jets of submanifolds}
\label{Jets of submanifolds}
 In view of the definition of the phase space, let us consider a
manifold
$\f M$
of dimension
$n$
and recall a few basic facts concerning jets of submanifolds
\cite{Vin88}.

 Let
$k \geq 0$
be an integer.
 A \emph{$k$--jet} of 1--dimensional submanifolds of
$\f M$
at
$x \in \f M$
is defined to be an equivalence class of 1--dimensional submanifolds
touching each other at
$x$
with a contact of order
$k \,.$
 The $k$--jet of a 1--dimensional submanifold
$s : \f N \hto \f M$
at
$x \in \f N$
is denoted by
$j_ks(x) \,.$
 The set of all $k$--jets of all 1-dimensional submanifolds at
$x \in \f M$
is denoted by
$J_{k \, x}(\f M,1) \,.$
 The set
$J_k(\f M,1) \byd \bigsqcup_{x \in \f M} J_{k \, x}(\f M,1)$
is said to be the \emph{$k$--jet space} of 1--dimensional
submanifolds of
$\f M \,.$
 In particular, for
$k = 0 \,,$
we have the natural identification
$J_0(\f M,1) = \f M \,,$
given by
$j_0s(x) = x \,,$
for each 1--dimensional submanifold
$s : \f N \hto \f M \,.$
 For each integers
$k \geq h \geq 0 \,,$
we have the natural projection
$\pi^k_h : J_k(\f M,1) \to J_h(\f M,1) :
j_ks(x) \mto j_hs(x) \,.$

 A chart of
$\f M$
is said to be \emph{divided} if the set of its coordinate functions
is divided into two subsets of 1 and
$n-1$
elements.
 Our typical notation for a divided chart will be
$(x^0,x^i) \,,$
with
$1 \le i \le n-1 \,.$
A divided chart and a 1--dimensional submanifold
$s : \f N \hto \f M$
are said to be \emph{related} if the map
$\br x^0 \byd x^0|_\f N \in \map(\f N, \, \Rn)$
is a chart of
$\f N \,.$
 In such a case, the submanifold
$\f N$
is locally characterised by
$s^i \com (\br x^0)^{-1} \byd
(x^i \com s) \com (\br x^0)^{-1} \in \map(\Rn, \Rn) \,.$
 In particular, if the divided chart is adapted to the submanifold,
then the chart and the submanifold are related.

 Let us consider a divided chart
$(x^0, x^i)$
of
$\f M \,.$

 Then, for each submanifold
$s : \f N \hto \f M$
which is related to this chart, the chart yields naturally the local
fibred chart
$(x^0, x^i; \, x^i_\ul \alp)_{1 \leq |\ul\alp| \leq k} \in
\map(J_k(\f M,1), \; \Rn^n \car \Rn^{k(n-1)})$
of
$J_k(\f M,1) \,,$
where
$\ul\alp \byd (h)$
is a multi--index of ``range" 1 and ``length"
$|\ul\alp| = h$
and the functions
$x^i_\ul\alp$
are defined by
$x^i_\ul\alp \com j_1\f N \byd \der_{0 \dots 0} \, s^i \,,$
with
$1 \leq |\ul\alp| \leq k \,.$

 We can prove the following facts:

1) the above charts
$(x^0, x^i; \, x^i_{\ul\alp})$
yield a smooth structure of
$J_k(\f M, 1) \,;$

2) for each 1--dimensional submanifold
$s : \f N \sub \f M$
and for each integer
$k \geq 0 \,,$
the subset
$j_ks(\f N) \sub J_k(\f M, 1)$
turns out to be a smooth 1--dimensional submanifold;

3) for each integers
$k \geq h \geq 1 \,,$
the maps
$\pi^k_h : J_k(\f M,1) \to J_h(\f M,1)$
turn out to be smooth bundles.

 We shall always refer to such divided charts
$(x^0, x^i)$
of
$\f M$
and to the induced fibred charts
$(x^0, x^i; \, x^i_\ul\alp)$
of
$J_k(\f M,1) \,.$

 Let
$m_1 \in J_1(\f M,1) \,,$
with
$m_0 = \pi^1_0 (m_1) \in \f M \,.$
Then, the tangent spaces at
$m_0$
of all 1--dimensional submanifolds
$s : \f N \hto \f M \,,$
such that
$j_1s(m_0) = m_1 \,,$
coincide.
Accordingly, we denote by
$T[m_1] \sub T_{m_0} \f M$
the tangent space at
$m_0$
of the above 1--dimensional submanifolds
$\f N$
generating
$m_1 \,.$
 We have the natural fibred isomorphism
$J_1(\f M,1) \to \Grass(\f M,1) :
m_1 \mto T[m_1]  \sub T_{m_0} \f M$
over
$\f M$
of the 1st jet bundle with the Grassmannian bundle of dimension 1.
 If
$s : \f N \hto \f M$
is a 1--dimensional submanifold, then we obtain
$T[j_1s] = \Span\lang \der_0 + \der_0s^i \der_i \rang \,,$
with reference to a related chart.
\subsubsection{Phase space}
\label{E: Phase space}
 A \emph{motion} is defined to be a 1--dimensional timelike
submanifold
$s : \f T \hto \f E \,.$

 For every arbitrary choice of a ``\emph{proper time origin}"
$t_0 \in \f T \,,$
we obtain the ``\emph{proper time scaled function}" given by the
equality
$\sig : \f T \to \baB T :
t \mto \fr1c \int_{[t_0, t]} \|\fr{ds}{d\br x^0}\| \, d\br x^0 \,.$

 This map yields, at least locally, a bijection
$\f T \to \baB T \,,$
hence a (local) affine structure of
$\f T$
associated with the vector space
$\baB T \,.$
 Indeed, this (local) affine structure does not depend on the choice
of the proper time origin and of the spacetime chart.

 Let us choose a time origin
$t_0 \in \f T$
and consider the associated proper time scaled function
$\sig : \f T \to \baB T$
and the induced linear isomorphism
$T\f T \to \f T \car \baB T \,.$
 Moreover, let us consider a spacetime chart
$(x^\lam)$
and the induced chart
$(\br x^0) \in \map(\f T, \Rn) \,.$
 Let us set
$\der_0 s^\lam \byd \fr{d s^\lam}{d\br x^0} \,.$

 The \emph{1st differential} of the motion
$s$
is defined to be the map
$ds \byd \fr{ds}{d\sig} : \f T \to \B T^* \ten T\f E \,.$

 We have
$g(ds, \, ds) = - c^2 \,.$

\medskip

 We assume as \emph{phase space} the subspace 
$\M J_1\f E \sub J_1(\f E,1)$
consisting of all 1--jets of motions.

 For each 1--dimensional submanifold
$s : \f T \sub \f E$
and for each
$x \in \f T \,,$
we have
$j_1s(x) \in \M J_1\f E$
if and only if
$T[j_1s(x)] = T_x\f T$
is timelike.

 Any spacetime chart
$(x^0, x^i)$
is related to each motion
$s \,.$
 Hence, the fibred chart
$(x^0, x^i, x^i_0)$
is defined on tubelike open subsets of
$\M J_1\f E \,.$

 We shall always refer to the above fibred charts.

\smallskip

 The \emph{velocity} of a motion
$s$
is defined to be its 1--jet
$j_1s : \f T \to \M J_1(\f E,1) \,.$

 We define the \emph{contact map} to be the unique fibred morphism
$\K d : \M J_1\f E \to \bar{\B T}^* \ten T\f E$
over
$\f E \,,$
such that
$\K d \com j_1s = ds \,,$
for each motion
$s \,.$
 We have
$g \, (\K d, \K d) = - c^2 \,.$
 The coordinate expression of
$\K d$
is
\beq
\K d =
c \, \alp^0 \,  (\der_0 + x^i_0 \, \der_i) \,,
\ssep{where}
\alp^0 \byd
1 /\sqrt{|g_{00} + 2 \, g_{0j} \, x^j_0 + g_{ij} \, x^i_0 \, x^j_0|}
\,.
\eeq

 The map
$\K d : \M J_1\f E \to \B T^*\ten T\f E$
is injective.
 Indeed, it makes
$\M J_1\f E \sub \B T^*\ten T\f E$
the fibred submanifold over
$\f E$
characterised by the constraint
$g_{\lam\mu} \, \dt x^\lam_0 \, \dt x^\mu_0 = - (c_0)^2 \,.$

\smallskip

 We define the \emph{time form} to be the map
$\tau
\byd - \fr1{c^2} g\Fla(\K d) :
\M J_1\f E \to \bar{\B T}\ten T^*\f E \,.$
 We have
$\tau (\K d) = 1$
and 
$\ba g(\tau,\tau) = - \fr1{c^2} \,.$
 The coordinate expression of
$\tau$
is
\beq
\tau = \tau_\lam \, d^\lam \,,
\ssep{where}
\tau_\lam = - \fr{\alp^0}{c}
\, (g_{0\lam}+ g_{i\lam}\, x^i_0) \,.
\eeq

\smallskip

 We define also the \emph{complementary contact map}
$\tht \byd 1 - \K d \ten \tau :
\M J_1\f E \to T^*\f E \ten T\f E \,.$
 The coordinate expression of
$\tht$
is
\beq
\tht =
d^\lam \ten \der_\lam +
(\alp^0)^2 \, (g_{0\lam} + g_{i\lam} \, x^i_0) \,
d^\lam \ten (\der_0 + x^j_0 \, \der_j) \,.
\eeq
\subsubsection{Vertical bundle of the phase space}
\label{E: Vertical bundle of the phase space}
 Let
$V\M J_1\f E \sub T\M J_1\f E$
be the vertical tangent subbundle over
$\f E \,.$
 The vertical prolongation of the contact map yields the mutually
inverse linear fibred isomorphisms
\beq
\nu_\tau : \M J_1\f E \to \B T \ten V^*_\tau\f E \ten V\M J_1\f E
\ssep{and}
\nu^{-1}_\tau :
\M J_1\f E \to V^*\M J_1\f E \ten \B T^* \ten V_\tau\f E \,,
\eeq
with coordinate expressions
\beq
\nu_\tau =  \fr1{c \, \alp^0} \, (d^i- x^i_0\, d^0) \ten \der^0_i\,,
\quad
\nu^{-1}_\tau = c \, \alp^0 \, d^i_0 \ten 
\big(\der_i - c\,\alp^0\tau_i(\der_0+x^p_0\,\der_p)\big) \,.
\eeq

 Thus, for each
$Y \in \sec(\M J_1\f E, V\M J_1\f E)$ 
and 
$X \in \sec(\f E,T\f E) \,,$ 
we obtain
\beq
\nu^{-1}_\tau(Y) \in \fib(\M J_1\f E, \, \B T^* \ten V_\tau\f E)
\ssep{and}
\nu_\tau(X) \in \sec(\M J_1\f E, \, \B T \ten V\M J_1\f E) \,,
\eeq  
with coordinate expressions
\beq
\nu^{-1}_\tau(Y) = c \, \alp^0 \, Y^i_0 \, \big(\der_i 
- c\,\alp^0\tau_i(\der_0+x^p_0\,\der_p)\big) 
\ssep{and}
\nu_\tau(X) = \fr 1{c\,\alp^0} \, \ti X^i \, \der^0_i \,,
\eeq
where $\ti X^i = X^i - x^i_0 \, X^0 \,.$
\subsubsection{Spacetime connections}
\label{E: Spacetime connections}
 We define a \emph{spacetime connection} to be a torsion free linear
connection 
$K : T\f E \to T^*\f E \ten TT\f E$
of the bundle
$T\f E \to \f E \,.$
 Its coordinate expression is of the type
\beq
K =
d^\lam \ten 
(\der_\lam + K\col\lam\nu\mu \, \dot x^\mu \, \dt\der_\nu) \,,
\ssep{with}
K\col\mu\nu\lam = K\col\lam\nu\mu \in \map(\f E, \, \Rn) \,.
\eeq

\smallskip

  We denote by 
$K[g]$ 
the \emph{Levi Civita connection}, i.e. the torsion free
linear spacetime connection such that 
$\nab g=0 \,.$
\subsubsection{Phase connections}
\label{E: Phase connections}
 We define a \emph{phase connection} to be a connection of the
bundle
$\M J_1\f E \to \f E \,.$

A phase connection can be represented, equivalently, by a tangent
valued form\linebreak
$\Gam : \M J_1\f E \to T^*\f E \ten T\M J_1\f E \,,$
which is projectable over
$\1 : \f E \to T^*\f E \ten T\f E \,,$
or by the complementary vertical valued form
$\nu[\Gam] :
\M J_1\f E \to T^*\M J_1\f E \ten V\M J_1\f E \,,$
or by the vector valued form
$\nu_\tau[\Gam] \byd \nu^{-1}_\tau \com \nu[\Gam] : \M J_1\f E \to
T^*\M J_1\f E \ten (\B T^* \ten V_\tau\f E) \,.$
 Their coordinate expressions are
\bgt
\Gam =
d^\lam \ten(\der_\lam + \Gam\Ga\lam i \, \der_i^0) \,,
\qquad
\nu[\Gam] =
(d^i_0 - \Gam\Ga\lam i \, d^\lam) \ten \der^0_i \,,
\\
\nu_\tau[\Gam] =
c \, \alp^0 \, (d^i_0 - \Gam\Ga\lam i \, d^\lam) \ten 
	\big(\der_i - c\,\alp^0\tau_i(\der_0+x^p_0\,\der_p)\big)\,,
\ssep{with}
\Gam\Ga\lam i \in \map(\M J_1\f E, \, \Rn) \,.
\end{gather*}

 We define the \emph{curvature} of a phase connection
$\Gam$
to be the vertical valued 2--form
\beq
R=R[\Gam] \byd - [\Gam, \, \Gam] :
\M J_1\f E \to
\Lam^2 T^*\f E \ten V\M J_1\f E \,,
\eeq
where
$[\,,]$
is the Fr\"olicher--Nijenhuis bracket.

 We can prove that there is a natural map 
$\chi : K \mto \Gam$
between linear spacetime connections 
$K$ 
and phase connections 
$\Gam \,,$
with coordinate expression
\beq
\Gam\Ga \lam i = 
K\col\lam i 0 + K\col\lam i p \ x^p_0 - x^i_0 \,
(K\col\lam 0 0 + K\col\lam 0 p \, x^p_0) \,.
\eeq
\subsubsection{Dynamical phase connection}
\label{E: Dynamical phase connection}
 The space of 2--jets of motions
$\M J_2\f E$
can be naturally regarded as the affine subbundle
$\M J_2\f E \sub \B T^* \ten T\M J_1\f E \,,$
which projects on
$\K d : \M J_1\f E \to \B T^* \ten T\f E \,.$

 A \emph{dynamical phase connection} is defined to be a 2nd--order
connection, i.e. a section
$\gam : \M J_1\f E \to \M J_2\f E \,,$
or, equivalently, a section
$\gam : \M J_1\f E \to \B T^* \ten T\M{J}_1\f E \,,$
which projects on
$\K d \,.$

 The coordinate expression of a dynamical phase connection is of the
type
\beq
\gam = c \, \alp^0 \,
(\der_0 + x^i_0 \, \der_i + \gam\ga i \, \der^0_i) \,,
\ssep{with}
\gam\ga i \in \map(\M J_1\f E, \, \Rn) \,.
\eeq

 If
$\gam$
is a dynamical phase connection, then we have
$\gam \con \tau = 1 \,.$

 The contact map
$\K d$
and a phase connection
$\Gam$
yield the section
$\gam \eqv \gam[\K d, \Gam] \byd 
\K d \con \Gam : \M J_1\f E \to \B T^* \ten T\M J_1\f E \,,$
which turns out to be a dynamical phase connection, with coordinate
expression
\beq
\gam\ga i = \Gam\Ga 0i + \Gam\Ga ji \, x^j_0 \,.
\eeq

 In particular, a linear spacetime connection
$K$
yields the dynamical phase connection
$\gam \byd \gam[\K d, K] \byd \K d \con \chi(K) \,,$
with coordinate expression
\begin{gather*}
\gam\ga i =
\\
=
K\col 0i0 +
K\col 0ih \, x^h_0 +
K\col hi0 \, x^h_0 +
K\col hik \, x^h_0 \, x^k_0 -
x^i_0 \, (
K\col 000  +
K\col 00h \, x^h_0 +
K\col h00 \,  x^h_0 +
K\col h0k  \, x^h_0 \, x^k_0
) \,.
\end{gather*}
\subsubsection{Phase 2--form and 2--vector}
\label{E: Phase 2--form and 2--vector}
 The metric
$g$
and a phase connection
$\Gam$
yield the scaled 2--form 
$\Ome \,,$ 
called \emph{(scaled) phase
2--form}, and  the scaled vertical 2--vector 
$\Lam \,,$ 
called \emph{(scaled) phase 2--vector}, 
\bat{3}
\Ome 
&\byd 
\Ome[g,\Gam] 
&&\byd
g \con \big(\nu_\tau[\Gam]\wed\tht\big) 
&&:
\M J_1\f E \to (\B T^* \ten \B L^2) \ten \Lam^2T^*\M J_1\f E \,,
\\
\Lam 
&\byd 
\Lam[g, \Gam] 
&&\byd
\ba g \con (\Gam \wed \nu_\tau) 
&&:
\M J_1\f E \to (\B T \ten \B L^{-2}) \ten \Lam^2T\M J_1\f E \,,
\end{alignat*}
with coordinate expressions
\beq
\Ome =
c \, \alp^0 \, (g_{i\mu} + c^2 \, \tau_i \, \tau_\mu) \,
(d^i_0 - \Gam\Ga \lam i \, d^\lam) \wed d^\mu\,,
\quad
\Lam =
\fr1{c \, \alp^0} \, (g^{j\lam} - x^j_0 \, g^{0\lam}) \,
(\der_\lam + \Gam\Ga \lam i \, \der^0_i) \wed \der^0_j \,.
\eeq

 We can easily see that  
$- c^2\,\tau \wed \Ome^3 \not\equiv 0$ 
and
$- \tfr1{c^2} \, \gam \wed \Lam^3 \not\equiv 0 \,.$ 

 There is a unique dynamical phase connection
$\gam \,,$
such that
$\gam \con \Ome[g, \Gam] = 0 \,.$
Namely,
$\gam = \gam[\K d, \Gam] \,.$

 In particular, a metric and time preserving spacetime connection
$K$
yields the (scaled) phase 2--form 
$\Ome[g,K]\byd \Ome[g,\chi(K)]$ 
and the (scaled) phase 2--vector
$\Lam[g,K]\byd \Lam[g,\chi(K)]$
with coordinate expressions
\bal
\Ome
&=
- c \, (g_{i\mu} + c^2 \, \tau_i \, \tau_\mu) \,
\big(d^i_0 -
(K\col \lam i0 + K\col \lam ij \, x^j_0 -
K\col \lam 00 \, x^i_0 - K\col \lam 0j \, x^i_0 \, x^j_0) \,
d^\lam\big) \wed d^\mu \,,
\\
\Lam 
&=
\fr1{c \, \alp^0} \, (g^{h\lam} - g^{0\lam} \, x^h_0) \,
\big(\der_\lam 
+ (K\col \lam i 0 + K\col \lam i j \, x^j_0
- K\col \lam 0 0 \, x^i_0 - K\col\lam 0 j \, x^i_0 \, x^j_0) \, 
\der^0_i\big) \wed \der^0_h \,.
\end{align*}
\subsubsection{Dynamical structures of the phase space}
\label{E: Dynamical structures of the phase space}
 Let us consider a phase connection
$\Gam$
and the induced phase objects
$\gam \byd \gam[\K d, \Gam] \,,$
$\Ome \byd \Ome[g, \Gam] \,,$ 
and
$\Lam \byd \Lam[g, \Gam] \,.$

\smallskip

 We define the Lie derivatives
\beq
L_\Gam\tau = (i_\Gam d - di_\Gam)\tau
\ssep{and}
L_R\tau = (i_R d + d i_R)\tau\,.
\eeq

 Then, the following results holds \cite{JanMod06}.

\bTh
The following assertions are equivalent.

\smallskip

{\rm(1)}
$L_{\nu_\tau(X)} \, L_\Gam \, \tau = 0 \,,$
$\Al \, X \in \sec(\f E,T\f E ) \,,$
and
$L_R \, \tau = 0 \,.$

\smallskip

{\rm (2)} 
$d\Ome = 0 \,,$ 
i.e. 
$(- c^2 \, \tau, \, \Ome)$ 
is a (scaled) almost--cosymplectic--contact pair.

\smallskip

{\rm (3)}
$[- \tfr1{c^2} \, \gam, \, \Lam] =  
\tfr1{c^2} \, \gam \wed \Lam\Sha(L_{\gam} \, \tau))$ 
and 
$[\Lam, \, \Lam] = 
2 \, \gam \wed(\Lam\Sha \ten \Lam\Sha)(d\tau)) \,,$ 
i.e. 
$(-\tfr1{c^2} \, \gam, \, \Lam, \, -c^2 \, \tau)$ 
is a (scaled regular) almost--coPoisson--Jacobi 3--plet.

 Moreover, the almost--cosymplectic--contact pair 
$( - c^2 \, \tau, \, \Ome)$ 
and the (regular) almost--coPoisson--Jacobi 3--plet 
$(- \tfr1{c^2} \, \gam, \, \Lam, \, -c^2 \, \tau)$
are mutually dual.
\hfill\ENDE
\eTh

\bLm
 We have 
$$
\Ome - c^2 \, L_\Gam \, \tau = - c^2\, d\tau \,.
\eqno{\ENDE}
$$
\eLm

\bTh
 The following assertions are equivalent.

\smallskip

{\rm (1)} 
$L_\Gam \, \tau = 0 \,.$

\smallskip

{\rm (2)} 
$\Ome = - c^2\, d\tau \,,$ 
i.e. 
$(-c^2\,\tau,\Ome)$ 
is a (scaled) contact pair.

\smallskip

{\rm (3)} 
$[- \tfr1{c^2} \, \gam, \Lam] = 0$ 
and 
$[\Lam,\Lam] = \tfr2{c^2} \, \gam \wed \Lam \,,$ 
i.e. 
$(-\tfr1{c^2} \, \gam, \Lam)$ 
is a (scaled regular) Jacobi pair.

\smallskip

 Moreover, the contact pair 
$(- c^2 \, \tau, \Ome)$ 
and the (regular) Jacobi pair 
$(-\tfr1{c^2} \, \gam,\Lam)$ 
are mutually dual.
\hfill\ENDE
\eTh

\medskip

 Next, let us consider a linear spacetime connection
$K$ 
and the induced phase objects 
$\Gam \byd \chi(K) \,,$ 
$\gam \byd \gam[\K d, \Gam] \,,$
$\Ome \byd \Ome[g, \Gam] \,,$ 
and
$\Lam \byd \Lam[g, \Gam] \,.$

\bTh\label{Th3.6}
 The following assertions are equivalent.

\smallskip

{\rm (1)}
$L_{\chi(K)} \, \tau = 0 \,.$
\bml
(2)
\quad
g(Z,Z)\,\big((\nab_Xg)(Y,Z) - (\nab_Yg)(X,Z) + g(T(X,Y),Z)\big)
\\
+ \tfr12 g(Z,X)(\nab_Yg)(Z,Z) - \tfr12 g(Z,Y)(\nab_X g)(Z,Z) = 0 \,,
\end{multline*}
for each 
$X,Y,Z \in \sec(\f E, T\f E) \,,$ 
where 
$T$ 
is the torsion of 
$K \,.$

\smallskip

{\rm (3)}
$\Ome = - c^2 \, d\tau \,,$ 
i.e. 
$(-c^2 \, \tau, \Ome)$ 
is a (scaled) contact pair.

\smallskip

{\rm (4)}
$[-\tfr1{c^2} \, \gam, \Lam] = 0$ 
and 
$[\Lam, \Lam] = \tfr2{c^2} \, \gam \wed \Lam \,,$ 
i.e. 
$(-\tfr1{c^2} \, \gam, \Lam)$ 
is a (scaled regular) Jacobi pair.

\smallskip

 Moreover, if the above conditions are fulfilled, then the contact
pair 
$(- c^2\,\tau, \Ome)$ 
and the (regular) Jacobi pair 
$(-\tfr1{c^2}\gam,\Lam)$ 
are mutually dual.
\hfill\ENDE
\eTh

\bCr
 Let 
$K$ 
be a torsion free spacetime connection. 
 If 
$\nab g$ 
and 
$g \ten \nab g$ 
are symmetric
(0,3) and (0,5) tensor fields, respectively, then  
$( -c^2 \, \tau \,, \, \Ome)$  
and 
$(-\tfr1{c^2} \, \gam \,,\, \Lam)$ 
are mutually dual contact and Jacobi pairs, respectively.
\hfill\ENDE
\eCr

 Eventually, let us consider the Levi Civita spacetime connection
$K[g]$ 
and the induced phase objects
$\Gam \eqv \Gam[g] \byd \chi(K) \,,$ 
$\gam \eqv \gam[\K d, g] \byd \gam[\K d, \Gam] \,,$
$\Ome \eqv \Ome[g] \byd \Ome[g, \Gam] \,,$ 
and
$\Lam[g] \byd \Lam[g, \Gam] \,.$

 Then, the equality 
$\nab g = 0$
and Theorem \ref{Th3.6} yield the following result.

\bTh
 We have: 

\smallskip

{\rm (1)} 
$\Ome= - c^2 \, d\tau \,,$ 
i.e. 
$(-c^2 \, \tau, \Ome)$
is a (scaled) contact pair.

\smallskip

{\rm (2)} 
$[-\tfr1{c^2} \, \gam, \Lam] = 0$ 
and 
$[\Lam, \Lam] = \tfr2{c^2} \, \gam \wed \Lam \,,$
i.e. 
$(-\tfr1{c^2} \, \gam, \Lam)$ 
is a (scaled regular) Jacobi pair.

\smallskip

 Moreover, the contact pair 
$( - c^2 \, \tau, \Ome)$ 
and
the (regular) Jacobi pair 
$(-\tfr1{c^2} \, \gam, \Lam)$ 
are mutually dual.
\hfill\ENDE
\eTh

\bRm
 If
$K$
is a spacetime connection, then the induced pairs 
$(- c^2 \, \tau , \, \Ome)$ 
and
$(-\tfr{1}{c^2}\gam, \, \Lam)$ 
are scaled. 

 On the other hand, some results of the general theory of geometrical
structures developed in the first two sections requires unscaled
pairs.

 Indeed, if we refer to a particle of mass
$m \in \B M \, $ 
and consider the universal scales
$\h \in \B T^{-1} \ten \B L^2 \ten \B M$
and
$c \in \B T^{-1} \ten \B L \,,$
then we obtain unscaled pairs in the following natural way.

 We have the unscaled spacetime 1--form
\beq
- \tfr{m\,c^2}{\h} \, \tau : \f E \to T^*\f E \,.
\eeq

 Moreover, the rescaled contact map
$\K D \byd \tfr{\h}{m\, c^2} \, \K d : \M J_1\f E \to T\f E$
yields the unscaled phase vector field
\beq
- \gam[\K D,K] = 
- \tfr{\h}{m\, c^2} \, \gam[\K d,K] : \f E \to T\M J_1\f E \,.
\eeq

 Furthermore, the rescaled metric
$G \byd \tfr m{\h} \, g : \f E \to \B T \ten T^*\f E \ten T^*\f E$
yields the unscaled phase 2--form and phase 2--vector
\bal
\Ome \eqv \Ome[G,K] 
&= 
\tfr m{\h} \, \Ome[g,K] : \M J_1\f E \to \Lam^2T^*\M J_1\f E \,,
\\
\Lam \eqv \Lam[G,K] 
&= 
\tfr \h{m} \, \Lam[g,K] : \M J_1\f E \to \Lam^2T\M J_1\f E \,.
\end{align*}

 Thus, if
$K$
is the Levi Civita spacetime connection, then
$(- \tfr{m\,c^2}{\h} \, \tau, \, \tfr{m}{\h} \, \Ome)$
and \linebreak
$(-\tfr{\h}{m\,c^2} \, \gam, \, \tfr{\h}{m} \, \Lam)$
turn out to be mutually dual unscaled contact and Jacobi
pairs of the phase space.

 Indeed, the Plank constant does not play any direct role in
classical mechanics; nevertheless, such a scale is necessary for
getting unscaled objects as above.
\hfill\ENDE
\eRm
{\small

}
\end{document}